\documentclass[12pt,oneside]{amsart}
\usepackage{geometry}                
\usepackage{graphicx}
\usepackage{amssymb}
\usepackage{epstopdf}

\usepackage{amsmath,amsthm,amscd,amssymb, enumerate}
\usepackage{latexsym}
\usepackage[colorlinks,citecolor=red,pagebackref,hypertexnames=false]{hyperref}
\numberwithin{equation}{section}
\theoremstyle{plain}
\newtheorem{theorem}{Theorem}[section]

\newtheorem{corollary}[theorem]{Corollary}
\newtheorem{proposition}[theorem]{Proposition}

\theoremstyle{definition}
\newtheorem{definition}[theorem]{Definition}

\theoremstyle{remark}
\newtheorem{remark}[theorem]{Remark}

\newtheorem{case[theorem]}{Case}

\def\bc{\begin{corollary}}
\def\ec{\end{corollary}}
\def\be{\begin{equation}}
\def\ee{\end{equation}}
\def\bast{\begin{eqnarray*} }
\def\east{\end{eqnarray*} }
\def\bea{\begin{eqnarray}}
\def\eea{\end{eqnarray}}

\def\ss{\smallskip}
\def\ms{\medskip}
\def\bs{\bigskip}

\def\ni{\noindent}

\def\R{\mathbb R}
\def\hd{{\dim_{\mathcal H}}}

\def\R{\Bbb R}

\def\bx{{\bf x}}

\def\vt{{\bf{t}}\,}

\def\dt{d_{\text{tot}}}
\def\dl{d_L^\sigma}
\def\dr{d_R^\sigma}
\def\meff{m_{\text{eff}}}
\def\pk{\mathcal P_k}

\def\hs{\hat{\sigma}}

\def\hd{\dim_{\mathcal H}}

\def\({\left(}
\def\){\right)}
\def\[{\left[}
\def\]{\right]}
\def\<{\left\langle}
\def\>{\right\rangle}

\def\zts{Z_t^\sigma}

\def\zvt{Z_{\vt}}
\def\zvts{\zvt^{\sigma}}
\def\rvt{\mathcal R_{\vt}}
\def\rvte{\mathcal R_{\vt,\epsilon}}
\def\rvts{\rvt^\sigma}

\def\sd{\mathbb S^{d-1}}

\def\sd{\mathbb S^{d-1}}

\def\z{\mathbf 0}
\def\mzs{\setminus\z}

\def\corank{\hbox{corank}}
\def\dim{\hbox{dim}}

\def\inter{\hbox{int}}
\def\rank{\hbox{rank}}

\title[On $k$-point configuration sets with nonempty interior]{On $k$-point configuration sets \\ with nonempty interior}

\date{Version of April 26, 2021. AG supported by US NSF DMS-1906186.}    

\author{Allan Greenleaf, Alex Iosevich and Krystal Taylor}

\email{allan@math.rochester.edu}

\email{iosevich@math.rochester.edu}

\email{taylor.2952@osu.edu}

\address{Department of Mathematics, University of Rochester, Rochester, NY 14627}

\address{Department of Mathematics, University of Rochester, Rochester, NY 14627}

\address{Department of Mathematics, Ohio State University, Columbus, Ohio 43210}

\begin{document}
\maketitle

\begin{abstract} We give conditions for  $k$-point configuration sets of thin sets to have nonempty interior,  
  applicable to a wide variety of configurations.
This is a continuation of our earlier work  \cite{GIT19}  on $2$-point configurations,
extending  a theorem of Mattila and Sj\"olin \cite{MS99}
for distance sets in Euclidean spaces. 
We show that for a general class of $k$-point configurations,
the configuration set of a $k$-tuple of sets, $E_1,\,\dots,\, E_k$,
has nonempty interior provided that the sum of their Hausdorff dimensions   satisfies a lower bound,
dictated by optimizing $L^2$-Sobolev  estimates of associated generalized Radon transforms
over all nontrivial partitions of the $k$ points into two  subsets.
We illustrate the general theorems with numerous specific examples.
Applications  to 3-point configurations  include areas of triangles   in $\R^2$ or the radii of their circumscribing circles;
volumes of pinned parallelepipeds in $\R^3$; and ratios of pinned distances in $\R^2$ and $\R^3$. 
Results for 4-point configurations  include cross-ratios on $\R$,   triangle area pairs determined by quadrilaterals in $\R^2$,
and dot products of differences in $\R^d$.
\end{abstract}  

\maketitle


\section{Introduction} \label{sec intro}

A classical result of Steinhaus \cite{Steinhaus} states that if $E\subset \R^d, d\ge 1,$ has positive Lebesgue measure, 
then the difference set $E-E\subset \R^d$ contains a neighborhood of the origin.
 $E-E$ can interpreted as the set of two-point configurations, $x-y$, of points of $E$ modulo the translation group.

Similarly,  in the context of the Falconer distance set problem,  
a  theorem of Mattila and Sj\"olin \cite{MS99} states that if $E\subset\R^d,\, d\ge 2$, is compact,
then the distance set of $E$, 
 $\Delta(E)=:\{\, |x-y|: x,\, y\in E\}\subset\R$,
 contains an open interval,  i.e., has nonempty interior,
if  the Hausdorff dimension $\hd(E)>\frac{d+1}2$. 
This represented a strengthening of Falconer's original result \cite{Falc86},
from $\Delta(E)$ merely having positive Lebesgue measure to having nonempty interior,  for the same range of $\hd(E)$.
 This was generalized   to distance sets  with respect to  norms on $\R^d$ 
 with positive curvature unit spheres in  Iosevich, Mourgoglou and Taylor \cite{IMT11}.

These latter types of results, for two-point configurations in \emph{thin sets},
 i.e.,   $E$ allowed to have Lebesgue measure zero 
but satisfying a lower bound on $\hd(E)$,
were  extended by the current authors to 
more general settings in \cite{GIT19}: 
(i)  configurations in $E$ as measured by a general class of $\Phi$\emph{-configurations}, 
which can be vector-valued and nontranslation-invariant;
and (ii)   \emph{asymmetric} configurations, i.e.,  between points in sets $E_1$ and $E_2$ lying in different spaces, for example
between points and   circles in $\R^2$, or points and hyperplanes in $\R^d$. 

 We point out that there are a number of other results that  are explicitly, 
 or can be interpreted as being, concerned with establishing conditions
 under which configuration sets of  thin sets,
 have nonempty interior, including \cite{BIT16,CLP16,GIP14,IT19} in the continuous setting and 
\cite{BIP17,ST17,PPSVV17,PPV19,B19} in finite field analogues. 
\medskip

The purpose of the current paper is to extend our results in \cite{GIT19} from 2-point to quite general $k$-point configuration sets,  for $k\ge 3$,
using that paper's  Fourier integral operator (FIO) approach, making use of linear $L^2$-Sobolev estimates,
but now optimizing over all possible nontrivial partitions of the $k$ points into two subsets. 
The FIO method we describe works 
in the absence of symmetry and on general manifolds, 
and indeed,   exploring that generality, 
rather than  sharpness of the lower bounds on the Hausdorff dimensions,
is the focus of  the current work.

However, we will start by illustrating the variety of what can be obtained via this approach 
with configurations defined by  classical geometric quantities
in low dimensional Euclidean spaces.
We describe a number of concrete examples,
but emphasize that the choice of these specific configurations is arbitrary;
our general results can be applied  to other configurations of interest, 
with the Hausdorff dimension threshold guaranteeing that the configuration set has nonempty interior depending on 
the outcome of optimizing over a family of FIO estimates; see Thms. \ref{thm threepoint} and \ref{thm kpoint} for exact statements.
Our first example is the following.
\ms

\begin{theorem}\label{thm areas} (Areas and circumradii of triangles) 
If $E\subset\R^2$ is compact with $\hd(E)>5/3$, then 
\ms

(i) the set  of areas of triangles determined by triples of points of $E$, 
\be\label{def areas}
\left\{\frac12\left|\det\left[x-z,\, y-z\right]\right|\, : x,y,z\in E\right\}\subset \R,
\ee
contains an open interval; and
\ms

(ii)  the set  of radii of circles determined by triples of points in $E$ 
contains an open interval.
\end{theorem}

We will see in Remark \ref{rem two d pinned} that the FIO method does not yield a pinned version of Theorem \ref{thm areas}, 
i.e., says nothing about  two-point configuration sets, 
$$\left\{\frac12\left|\det\left[x-z,\, y-z\right]\right|\, : x,y\in E\right\},$$ 
for a fixed $z\in E$. 
However,   in dimension $d\ge 3$ it does yield a result for {\it all} $z\in E$ (or, indeed, any $z\in\R^d$),
and we have the following $d$-point configuration result. 
\ms

\begin{theorem}\label{thm volumes} (Strongly pinned\footnote{The term {\it pinned} is often used to refer to estimates for the supremum over $x^0\in E$ of expressions such as \eqref{def vols}. Here we obtain a result valid for {\it all} $x^0$, hence our adoption of {\it strongly pinned} for lack of a better term.} volumes) 
Let $d\ge 3$. If $E\subset\R^d$ is compact, then for any $x^0\in\R^d$, the set of  volumes of parallelepipeds determined by  $x^0$ and
$d$-tuples of points of  $E$,
\be\label{def vols}
V_d^{x^0}(E):=\left\{\left|\det\left[x^1-x^0, x^2-x^0,\dots, x^d-x^0\right]\right|\, : x^1,x^2,\dots, x^d\in E\right\},
\ee
has nonempty interior in $\R$ if $\hd(E)>d-1+(1/d)$.
\end{theorem}

\begin{remark} For $d=3$, this improves upon an earlier result of the first two authors and Mourgoglou \cite{GIM15}, 
which was that if $\hd(E)>13/5$, 
then $V_3^0(E)$ has positive Lebesgue measure.
\end{remark}
\ms

Returning to three-point configurations, our method also yields a result about ratios of distances.
S. Mkrtchyan and the first two named authors  studied in \cite{GIMk18}  the existence of similarities
of $k$-point configurations in thin sets. They posed the question of whether, under some lower bound  restriction 
on $\hd(E)$, for every $r>0$ there exist $x,y,z\in E$ such that $|x-z|=r|y-z|$. 
We can partially address this, showing that the set of such $r$ at least contains an interval.
\medskip

To put this in perspective, note that an immediate consequence of the result of Mattila and Sj\"olin \cite{MS99} is that
if $\hd(E)>(d+1)/2$, then
\be\label{eqn sjolin}
\inter\left(\left\{\frac{|w-z|}{|x-y|}\, :\, x,y,z,w\in E\right\}\right)\ne\emptyset.
\ee
(See also \cite{IKP18} for a finite field analogue.)
On the other hand, Peres and Schlag \cite{PSch00} showed that if $\hd(E)>(d+2)/2$, a stronger
property holds: 
\be\label{eqn schlag}
\hbox{ there exists an } x\in E\, \hbox{ s.t. }\, \inter\left(\left\{\frac{|x-z|}{|x-y|}\, :\, y,z\in E\right\}\right)\ne\emptyset.
\ee
(See also \cite{ITU16,IL19} for extensions of this.)
Here we prove a result for a property of intermediate strength, one which implies \eqref{eqn sjolin} but is in turn implied by \eqref{eqn schlag};
however,  in dimensions $d=2,3$ our result is proved for lower $\hd(E)$ than the known range for \eqref{eqn schlag}:
\ms

\begin{theorem}\label{thm distsim} 
(Ratios of pairs of  pinned distances) Let $d\ge 2$ and $E\subset \R^d$ compact.  Then, if $\hd(E)>(2d+1)/3$,

\be\label{eqn distsim}
\inter\left(\left\{\frac{|x-z|}{|x-y|}\, :\, x,y,z\in E,\, x\ne y \right\}\right)\ne\emptyset.
\ee
\end{theorem}
\bs

We now turn from three-point configurations to a pair of results concerning  four-point configurations, in $\R$ and $\R^2$, resp.

\begin{theorem}\label{thm cross} (Cross ratios) Let $E\subset \R$ be compact with $\hd(E)>3/4$.
Then the set of cross ratios of four-tuples of points of $E$,
$$Cross(E)=\left\{\left[x_1,x_2;x_3,x_4\right]=\frac{(x_3-x_1)(x_4-x_2)}{(x_3-x_2)(x_4-x_1)}\, :\, x_1,x_2,x_3,x_4\in E\right\}\subset\R,$$
contains an open interval.
\end{theorem}

So far, all of the configurations described have been measured by scalar-valued functions. 
Returning to $d=2$, an example of a {\it vector-valued} configuration is a variation of Theorem \ref{thm areas}, 
where one takes four points in the plane, say $x,y,z,w$, and considers the quadrilateral they generate.
Pick one of the two diagonals, say $\overline{yw}$; this  splits the quadrilateral into two triangles,
and we study the vector-valued configuration consisting of their areas.

\begin{theorem}\label{thm two areas} (Pairs of areas of triangles)
If $E\subset\R^2$ is a compact set  and  has $\hd(E)>7/4$, then the set of pairs of areas of triangles determined by 4-tuples of points of $E$,
\be\label{def two areas}
\left\{\left( \frac12\left|\det\left[x-w,y-w\right]\right|,\, \frac12\left|\det\left[ y-w,z-w\right]\right|\right):\, x,y,z,w\in E \, \right\}.
\ee
has nonempty interior in $\R^2$.
\end{theorem}
\bs

Finally, we give two more applications of the FIO method, this time  to configurations with a more additive combinatorics flavor. 
There has been considerable work on products of differences in the discrete or finite field setting; just a few references are 
 \cite{HIS07,B13,BHIPR17,MurP18,MurPRRS19}. 
 An analogue of some of these results in the continuous setting is the following
 \ms
 
 \begin{theorem}\label{thm prod of diff}
 For $d\ge 1$ and $E\subset\R^d$ compact, 
the set of dot products of differences of points in $E$,
  $$\left\{(x-y)\cdot(z-w)\, :\, x,y,z,w\in E\, \right\}\subset \R,$$
 \ms
 
\ni  has nonempty interior if $\hd(E)>\frac{d}2+\frac14$.
 \end{theorem}
 
 Another result of sum-product type is 
 
\begin{theorem}\label{thm bilinear forms} (Generalized  sum-product sets) Let $Q_1,\dots,Q_l$ 
be  nondegenerate, symmetric bilinear forms  
on $\R^{d}$.
Suppose that $E_i\subset \R^{d}$ are  compact sets with $\hd(E_i)>\frac{d}2+\frac1{2l}$ for all $1\le i\le 2l$. Then the set of values
\be\label{def form sum set}
\Sigma_{\vec{Q}}\left(E_1,\dots,E_{2l}\right):=\left\{ \sum_{j=1}^l Q_j\left(x^{2j-1},x^{2j}\right): x^i\in E_i,\, 1\le i\le 2l\right\}\subset\R
\ee
has nonempty interior. In particular, taking all of the $Q_j(x,y)=x\cdot y$, under the same conditions on the  $\hd(E_i)$,  
the sum-(Euclidean  inner) product set of the $E_i$,
$$
\left\{\left(x^1\cdot x^2\right)+\left(x^3\cdot x^4\right)
+\cdots +\left(x^{2l-1}\cdot x^{2l}\right): x^i\in E_i,\, 1\le i\le 2l\right\},$$
has nonempty interior.
\end{theorem}

\begin{remark} This follows from a more general result allowing the forms to be on spaces of different dimensions; 
see Thm. \ref{thm general bilinear forms}.
\end{remark}


\section{Three-point configurations}\label{sec three point}

To describe a general class of $k$-point configurations which includes the examples above, 
we start by recalling the framework of $\Phi$-configuration sets,
introduced  by Grafakos, Palsson and the first two authors in  \cite{GGIP12}.
We used this approach to establish nonempty interior results for 2-point 
configuration sets in  the current article's prequel, \cite{GIT19}. 
To  minimize the notation, we initially describe these for 3-point configurations,  introducing the basic method and results, which will be extended to 
higher $k$ in Sec. \ref{sec bigger than three}.
\smallskip

A  {\it 3-point configuration function} is initially  a smooth $\Phi:\R^d\times\R^d\times\R^d\to \R^p$ (with $p\le\nolinebreak d$);
we use the notation $\Phi(x^1,x^2,x^3),\, x^j\in\R^d,\, j=1,2,3$.  
Since, for many problems of interest there are points,
often corresponding to degenerate configurations, 
where $\Phi$ has critical points or fails to be smooth,  
it is useful to restrict the domain  of $\Phi$.
Anticipating the extension to $k$-point configurations later, we label the three copies of $\R^d$  (or open subsets of $\R^d$) as
$X^1,\, X^2,\, X^3$.
Furthermore, for some applications it is useful to allow
the $X^j$ to be manifolds,  of possibly different dimensions, $d_j,\, j=1,2,3$.
Thus, in general we define a 3-point configuration function to be a mapping $\Phi:X^1\times X^2 \times X^3 \to T$,
where $T\subset\subset\R^p$, or even a \linebreak$p$-dimensional manifold,
containing the range of $\Phi$ on the compact sets of interest.
Function spaces on the $X^j$ are with respect to smooth densities,
which  do not play a significant role and therefore  are suppressed  in the notation.
\medskip

For compact sets $E_j\subset X^j,\, j=1,2,3$,  define
the 3-point {\it $\Phi$-configuration set} of $E_1,E_2,E_3$,
\be\label{def three point}
\Delta_\Phi\left(E_1,E_2,E_3\right):=\left\{\Phi\left(x^1,x^2,x^3\right): \, x^j\in E_j, \, j=1,2,3\right\}\subset T.
\ee
The goal is to find conditions on $\hd(E_j)$ ensuring that $\inter\left(\Delta_\Phi\left(E_1,E_2,E_3\right)\right)\ne\emptyset$.
\medskip

If the full differential $D_{x^1,x^2,x^3}\Phi$ has maximal rank ($=p$) everywhere, i.e.,
$\Phi$ is a submersion, then $\Phi$ is a  defining function for a family of smooth  surfaces in \linebreak$X:=X^1\times X^2\times X^3$, and  
for each $\vt\in T$, the level set
\be\label{def zt}
Z_{\vt}:=\big\{(x^1,x^2,x^3)\in X: \Phi(x^1,x^2,x^3)=\vt\big\}
\ee
is smooth and of codimension $p$ in $X$, and $\zvt$ depends smoothly on $\vt$.
(For $p=1$,  we denote $\vt$  by simply $t$.)
\smallskip

If $s_j<\hd(E_j)$, let $\mu_j$ be a Frostman measure on $E_j$ with finite $s_j$-energy.
(See the discussion in Sec. \ref{sec background} below.)
The choice of the $\mu_j$ induces a {\it configuration measure}, $\nu$,  on $T$, having various equivalent definitions, e.g.,
for $g\in C_0(T)$,
\be\label{def config meas}
\int_{T} g(\vt) \, d\nu(\vt) = \int\int\int_{E_1\times E_2\times E_3} g(\Phi(x^1,x^2,x^3))\, d\mu_1(x^1)\, d\mu_2(x^2)\, d\mu_3(x^3).
\ee
If one can show that $\nu$ is absolutely continuous with respect to Lebesgue measure, $d\vt$;  its density function is continuous; 
and  $\Delta_\Phi(E_1,E_2,E_3)$ 
is nonempty, then it follows that $\inter(\Delta_\Phi(E_1,E_2,E_3))\ne\emptyset$.
\ss

Following the general approach of \cite{GIT19}, but now exploiting the fact that we are studying 3-point, rather than 2-point, configurations,
we will derive the continuity of the density of $d\nu$ (denoted $\nu(\vt)$) from $L^2$-Sobolev mapping properties of any of three different  families of 
generalized Radon transforms associated to $\zvt$, as follows.
\ss

Write a nontrivial partition  of $\{1,2,3\}$ as $\sigma=(\sigma_L|\sigma_R)$,
grouping the variable(s) $x^i$ corresponding to  $i\in \sigma_L$  on the left and  the variable(s) corresponding to $i\in\sigma_R$ on the right. 
Due to the symmetry of  $L^2$-Sobolev estimates for FIOs  under adjoints, we may assume that $|\sigma_L|=2,\, |\sigma_R|=1$; 
furthermore, permutation within $\sigma_L$ is irrelevant, so up to  interchange of those two indices, 
there are three such partitions, $\sigma=(12|3),\, (13|2)$ and $(23|1)$. 
Corresponding to each of these, for each $\vt\in T$, partitioning and permuting the variables according to $\sigma$,  the surface $Z_{\vt}$  defines  incidence relations,
\bea\label{def three zt}
\zvt^{(12|3)}&:=&\{(x^1,x^2;x^3): (x^1,x^2,x^3)\in\zvt\}\subset (X^1\times X^2)\times X^3,\nonumber\\
\zvt^{(13|2)}&:=&\{(x^1,x^3;x^2): (x^1,x^3,x^2)\in\zvt\}\subset (X^1\times X^3)\times X^2,\\
\zvt^{(23|1)}&:=&\{(x^2,x^3;x^1): (x^2,x^3,x^1)\in\zvt\}\subset (X^2\times X^3)\times X^1.\nonumber
\eea
Each   $\zvt^{\sigma}$, with $\sigma=(ij|k)$,  defines an incidence relation from $X^k$ to $X^i\times X^j$,
and to this is associated a generalized Radon transform, 
$\rvt^{\sigma}$; all the $\rvt^{\sigma}$ have the `same' Schwartz kernel, namely the singular measure supported on $\zvt$,
$$\lambda_\vt:=\chi(x^1,x^2,x^3)\cdot \delta\left(\Phi\left(x^1,x^2,x^3\right)-\vt\right),$$
 except  that the order and grouping of the variables are dictated by $\sigma$. 
  I.e., the kernel of $\rvt^{(ij|k)}$ is $K_\vt^{(ij|k)}(x^i,x^j,x^k):=\lambda(x^1,x^2,x^3)$.
 (Here $\chi$ is a fixed cutoff function $\equiv 1$ on $E_1\times E_2\times E_3$ which plays no further role.)
 \ss

For each $\sigma$, we  can formulate the double fibration  condition, $(DF)_\sigma$,
standard in the theory of generalized Radon transforms and originating in the works of Gelfand;  Helgason \cite{Helgason};
and Guillemin and Sternberg \cite{GuSt,GuSt79},
namely that  the two spatial projections from $\zvt^\sigma$ have maximal rank,
namely 
 \be\label{def DF}
 (DF)_{\sigma}\qquad \pi_{ij}:\zvt^{\sigma}\to X^i\times X^j\hbox{ and }\pi_k:\zvt^{\sigma}\to X^k\hbox{ are submersions.}
 \ee
 This implies that not only does $\rvt^\sigma: \mathcal D(X^k)\to \mathcal E(X^i\times X^j)$, but also
$$\rvt^\sigma:\mathcal E'(X^k)\to \mathcal D'(X^i\times X^j),$$
defined weakly by
 $$\rvt^{\sigma}f(x^i,x^j)=\int_{\{x^k:\, \Phi(x^1,x^2,x^3)=\vt\}} f(x^k),$$
where the integral is with respect to the surface measure induced by $\lambda_\vt$ on the codimension $p$ surface  
$\left\{x^k:\, \Phi\left(x^1,x^2,x^3\right)=\vt\right\}\subset X^k$.
\ss

An alternate description of the
configuration measure defined by \eqref{def config meas} is in terms of the $\rvt^\sigma$;
this was stated and proved in  the case of 2-point configuration measures in \cite[Sec. 3]{GIT19}.
However, the proof there goes over with minor modifications to the case of $k$-point configurations,
and for completeness we give the argument for $k=3$ in Sec. \ref{subsec density} below.
Namely, as long as  the terms in the two arguments of the $\<\cdot,\cdot\>$ pairing below belong to Sobolev spaces 
 on which the bilinear pairing is continuous, $\nu$ has a  density given by
 \be\label{eqn pairing}
 \nu(\vt)=\<\rvt^{(ij|k)}(\mu_k),\mu_i\times \mu_j\>.
 \ee

Now, under the double fibration condition $(DF)_\sigma$, the generalized Radon transform $\rvt^{\sigma}$ is a Fourier integral operator 
(FIO) associated with a canonical relation
  $$C_{\vt}^{\sigma}\subset \left(T^*\left(X^i\times X^j\right)\setminus 0\right)\times\left(T^*X^k\setminus 0\right),$$
  where $C_{\vt}^\sigma=(N^*\zvts)'$, the (twisted) conormal bundle of $\zvts$ (see Sec. \ref{sec background}).
  All three of $\rvt^{(12|3)},\, \rvt^{(13|2)},\, \rvt^{(23|1)}$ are Fourier integral operators of the same order,
 $$m=0+\frac{1}2 p-\frac14(d_1+d_2+d_3)=\frac{p}2-\frac14 d^{tot},\quad d^{tot}:=d_1+d_2+d_3.$$
However, due to the (possibly) different dimensions, in order to understand the optimal estimates for the operators $\rvt^{\sigma}$,
one knows from standard FIO theory that   the estimates are conveniently  expressed in
terms of what we will call their {\it effective orders}, $\meff^\sigma$. These are defined by writing  $m$  in three different ways, accounting 
for the dimension differences $\left| \dim\(X^i\times X^j\)-\dim\(X^k\)\right|$:
\bast
m&=&\meff^{(12|3)}-\frac14 |d_1+d_2-d_3|,\\
m&=&\meff^{(13|2)}-\frac14 |d_1+d_3-d_2|,\hbox{ or }\\
m&=&\meff^{(23|1)}-\frac14|d_2+d_3-d_1|.
\east
In terms of the $\meff^\sigma$, the mapping properties of the operators $\rvt^\sigma$ can be described as
$$\rvt^\sigma: L^2_r \to L^2_{r-\meff^\sigma-\beta_{\vt}^\sigma},\, \forall r\in\R,$$
for  certain (possible) losses $\beta_{\vt}^\sigma\ge 0$. If, for some value $\vt_0\in T$,  $C_{\vt_0}^\sigma$ is {\it nondegenerate}, 
i.e., one of its two natural  projections to the left or right, 
$\pi_L$ or $\pi_R$, is of maximal rank (which implies that the other is as well), 
then $\beta_{\vt_0}^\sigma=0$, and  by structural stability of submersions this is also true for all $\vt$ near $\vt_0$
(see Sec. \ref{sec background}.)
Our basic  assumption is  that,  for at least one $\sigma$, there  is a known $\beta^\sigma\ge 0$ such that 
$\rvt^\sigma: L^2_r \to L^2_{r-\meff^\sigma-\beta^\sigma}$ uniformly for $\vt\in T$.
\smallskip

To simplify the arithmetic, assume that for all the $\sigma=(ij|k)$, we have $d_i+d_j\ge d_k$,
which includes the equidimensional case, $d_1=d_2=d_3$.
Then ${\meff}^{(ij|k)}=(p-d_k)/2$, and thus
our basic boundedness assumption is that, for at least one of the $\sigma$, 
\be\label{eqn basic three}
\rvt^\sigma: L^2_r(X^k) \to L^2_{r+\frac12(d_k-p)-\beta_{\vt}^\sigma}(X^i\times X^j)\hbox{ uniformly in }\vt.
\ee

At the start of the argument, the $s_j,\, j=1,2,3,$ were chosen to be any values such that $\hd(E_j)>s_j$, and
each $\mu_j$ has finite $s_j$ energy, so that  $\mu_j\in L^2_{(s_j-d_j)/2}(X^j)$.
An easy calculation with Sobolev norms shows that if  $u_j\in L^2_{r_j}(\R^{d_j})$ with $r_j\le 0,\, j=1,2,$ 
then $u_1\otimes u_2\in L^2_{r_1+r_2}(\R^{d_1+d_2})$, and this extends to compactly supported distributions on manifolds
(see Prop. \ref{prop Sob} below).
Thus, 
$$\mu_1\times \mu_2\in L^2_{(s_1+s_2-d_1-d_2)/2},\, \mu_1\times \mu_3\in L^2_{(s_1+s_3-d_1-d_3)/2},\hbox{ and }
\mu_2\times \mu_3\in L^2_{(s_2+s_3-d_2-d_3)/2}.$$

Combining all of  these considerations, and focusing on $\sigma=(12|3)$ for the moment, 
we see that the bilinear pairing in the expression \eqref{eqn pairing} for $\nu(\vt)$ is continuous if

$$\left(s_3-d_3\right)/2+(d_3-p)/2-\beta^{(12|3)}+\left(s_1+s_2-d_1-d_2\right)/2\ge 0,$$
i.e., 
$$s_1+s_2+s_3\ge d_1+d_2+p+2\beta^{(12|3)}.$$
The analogous calculation holds for whichever of the $\rvt^\sigma$ one knows estimates for, 
and the minimum over $\sigma$  of the 
right hand sides gives a 
sufficient condition for $\nu(\vt)$ to be continuous.
Thus, the set where $\nu(\vt)>0$ is an open set;
to conclude that  $\inter\left(\Delta_\Phi\left(E_1,E_2,E_3\right)\right)\ne\emptyset$, 
it suffices to show that
$\Delta_\Phi\left(E_1,E_2,E_3\right)$ itself is nonempty. 

As in \cite{GIT19}, this follows
by noting that what we have done above already implies the Falconer-type conclusion that $\Delta_\Phi\left(E_1,E_2,E_3\right)\subset\R^p$
has positive Lebesgue measure.
In fact, if $\big\{B(\vt_j,\epsilon_j)\big\}$ is any cover of $\Delta_\Phi\left(E_1,E_2,E_3\right)$,
one has
\bea
1&=&\mu_1(E_1)\cdot\mu_2(E_2)\cdot\mu_3(E_3)
=\left(\mu_1\times\mu_2\times\mu_3\right)\left(E_1\times E_2\times E_3\right)\nonumber \\
& & \nonumber \\
&\le&\left(\mu_1\times\mu_2\times\mu_3\right)\left(\Phi^{-1}\left(\bigcup_j B\left(\vt_j,\epsilon_j\right) \right)\right) \nonumber \\
& & \nonumber \\
&\le&  \sum_j \left(\mu_1\times\mu_2\times\mu_3\right)\left(\Phi^{-1}\left(B\left(\vt_j,\epsilon_j\right)\right)\right)\nonumber\\
&=&\sum_j \nu\left(B\left(\vt_j,\epsilon_j\right)\right)
\le C_\Phi\sum_j\epsilon_j^p\label{eqn Leb}
\eea
by \eqref{eqn epsilon est} below, so that $\sum_j |B(\vt_j,\epsilon_j)|_p\ge C_\Phi'$ is bounded below.
Hence $\Delta_\Phi(E_1,E_2,E_3)$  has positive $p$-dimensional Lebesgue measure and is therefore nonempty;
 by the continuity of $\nu(\vt)$, it in fact has nonempty interior. 
\medskip

Summarizing, we have established the following method for proving that 3-point configuration sets have nonempty interior:

\begin{theorem} \label{thm threepoint}
(i) With the notation and assumptions as above, define
$$s_\Phi=p+\min\left(d_1+d_2+2\beta^{(12|3)},\, d_1+d_3 +2\beta^{(13|2)},\, d_2+d_3+2\beta^{(23|1)}\right),$$
where the $\min$ is taken over those of the partitions $\sigma=(ij|k)$ for which 
\ms

(a) the double fibration condition $(DF)_\sigma$  
\eqref{def DF} holds,
and 
\ms

(b) one has uniform boundedness of the generalized Radon transforms $\rvt^\sigma$ with loss of $\,\le \beta^\sigma$ derivatives 
\eqref{eqn basic three}.
\ms

Then, if $E_j\subset X^j$ are compact sets with $\hd(E_1)+\hd(E_2)+\hd(E_3)>s_\Phi$, 
it follows that $\inter\left(\Delta_\Phi\left(E_1,E_2,E_3\right)\right)\ne\emptyset$.
\bs

(ii) In particular, suppose that $X^1=X^2=X^3=X$, with $dim(X)=d$, and 
there is a partition $\sigma=(ij|k)$ such that (a) holds and the canonical relations $C_{\vt}^\sigma$ are 
nondegenerate (so that $\beta^\sigma=0$). It follows that, if $E\subset X$ is compact with $\hd(E)>(2d +p)/3$, then
 $\inter\left(\Delta_\Phi\left(E,E,E\right)\right)\ne\emptyset$.
\end{theorem}
\smallskip


\section{Background material}\label{sec background}

We give a brief survey of the relevant facts needed in the paper, referring for more background and further details 
to  H\"ormander \cite{Hor71,Hor85} for  Fourier integral operator theory,
Mattila \cite{Mat95,Mat15} for geometric measure theory,
and \cite{GIT19} for the case of 2-point configurations. 
\ms

\subsection{Fourier integral operators}\label{subsec FIO}

Let $X$ and $Y$ be smooth manifolds of dimensions $n_1,\, n_2$, resp. Then $T^*X,\, T^*Y$ are each symplectic manifolds, 
with canonical two-forms denoted $\omega_{T^*X},\, \omega_{T^*Y}$, resp. Equip $T^*X\times T^*Y$ 
with the {\it difference symplectic form}, $\omega_{T^*X} - \omega_{T^*Y}$. For our purposes, a {\it canonical relation} will mean a  submanifold, 
$C\subset (T^*X\mzs)\times(T^*Y\mzs)$  (hence of dimension $n_1+n_2$), which is  conic Lagrangian with respect to $\omega_{T^*X} - \omega_{T^*Y}$.
\ms

For some $N\ge 1$, let $\phi: X\times Y\times (\R^N\mzs)\to \R$ be a smooth phase function 
which is positively homogeneous of degree 1 in $\theta\in\R^N$, i.e., 
$\phi(x,y,\tau\theta)=\tau\cdot\phi(x,y,\theta)$ for all $\tau\in\R_+$. 
Let $\Sigma_\phi$ be the {\it critical set} of $\phi$ in the $\theta$ variables,
$$\Sigma_\phi:=\{(x,y,\theta)\in X\times Y\times (\R^N\mzs): d_\theta\phi(x,y,\theta)=0\},$$
and 
$$C_\phi:=\{(x,d_x\phi(x,y,\theta); y, -d_y\phi(x,y,\theta)): (x,y,\theta)\in\Sigma_\phi\},$$
both of which are conic sets.
If we impose  the first order nondegeneracy conditions
$$d_{x}\phi(x,y,\theta)\ne 0\hbox{ and } d_{y}\phi(x,y,\theta)\ne 0, \forall (x,y,\theta)\in \Sigma_\phi,$$
then $C_\phi\subset (T^*X\mzs)\times(T^*Y\mzs)$.
If in addition one demands that
$$\rank[d_{x,y,\theta}d_\theta\phi(x,y,\theta)]=N,\, \forall\, (x,y,\theta)\in\Sigma_\phi,$$
then $\Sigma_\phi$ is smooth, $dim(\Sigma_\phi)=n_1+n_2$, and the map
\be\label{eqn cphi param}
\Sigma_\phi\owns (x,y,\theta)\to \left(x,d_x\phi\left(x,y,\theta\right); y,-d_y\phi\left(x,y,\theta\right)\right)\in C_\phi
\ee
is an immersion, whose image is an immersed canonical relation; the phase function $\phi$ is said to {\it parametrize} $C_\phi$.
\ss

For a canonical relation $C\subset (T^*X\mzs)\times(T^*Y\mzs)$  and $m\in\R$,
one  defines $I^m(X,Y;C)=I^m(C)$, the class of  {\it Fourier integral operators} $A:\mathcal E'(Y)\to \mathcal D'(X)$ of order \nolinebreak$m$,
as the collection of operators whose Schwartz kernels are locally finite sums of oscillatory integrals of the form
$$K(x,y)=\int_{\R^N} e^{i\phi(x,y,\theta)} a(x,y,\theta)\, d\theta,$$
where $a(x,y,\theta)$ is a symbol of order $m-N/2+(n_1+n_2)/4$ and $\phi$ is a phase function as above, 
parametrizing some relatively open $C_\phi\subset C$.
\ms

The  FIO relevant for this paper are the {\it generalized Radon transforms} $\rvt$ 
determined by  defining functions $\Phi:X\times Y\to \R^p$ 
satisfying the double fibration condition that $D_x\Phi$ and $D_y\Phi$ have maximal rank.   The Schwartz kernel of each $\rvt$ is a smooth multiple of $\delta_p(\Phi(x,y)-\vt)$, 
where $\delta_p$ is the delta distribution on $\R^p$. From the Fourier inversion representation of $\delta_p$, we see that $\rvt$ has kernel
$$K_{\vt}(x,y)=\int_{\R^k} e^{i(\Phi(x,y)-\vt)\cdot\theta} \, b(x,y)\cdot 1(\theta)\, d\theta,$$
where $b\in C_0^\infty$. Since the amplitude is a symbol of order 0, $\rvt$ is an FIO of order $0+p/2-(n_1+n_2)/4=-(n_1+n_2-2p)/4$
associated with the canonical relation parametrized as in \eqref{eqn cphi param}
by $\phi(x,y,\theta)=(\Phi(x,y)-\vt)\cdot\theta$, which is the twisted conormal bundle
of the incidence relation $\zvt$,
$$C_{\vt}=N^*\zvt':=\left\{\left(x,\sum_{j=1}^k d_x\Phi_j\left(x,y\right)\theta_j; 
y,-\sum_{j=1}^k d_x\Phi_j\left(x,y\right)\theta_j\right): \left(x,y\right)\in \zvt,\, \theta\in\R^k\mzs\right\}.$$
For $T$-valued defining functions $\Phi$,
 as in the general formulation of our results, 
 this discussion is easily modified by introducing local coordinates on $T$.
\ss

For a general canonical relation, $C$, the natural projections $\pi_L:T^*X\times T^*Y \to T^*X$ and $\pi_R:T^*X\times T^*Y \to T^*Y$ restrict to 
$C$, and by abuse of notation we refer to the restricted maps with the same notation.
One can show that, at any point $c_0=(x_0,\xi_0;y_0,\eta_0)\in C$, one has $\corank(D\pi_L)(c_0)=\corank(D\pi_R)(c_0)$; we say that
the canonical relation $C$ is {\it nondegenerate} if this corank is zero at all points of $C$, i.e., if
$D\pi_L$ and $D\pi_R$ are of maximal rank. 
If $\dim(X)=\dim(Y)$, then $C$ is nondegenerate iff $\pi_L,\, \pi_R$ are local diffeomorphisms, and then $C$ is a {\it local canonical graph}, i.e., locally 
near any $c_0\in C$ is equal to the graph of a canonical transformation. If $\dim(X)=n_1>n_2=\dim(Y)$, then $C$ is nondegenerate iff $\pi_L$ is an 
immersion and $\pi_R$ is a submersion.
To describe the $L^2$-Sobolev estimates for FIOs, it is convenient to normalize the order and 
consider $A\in I^{\meff-\frac{|m_1-m_2|}4}(C)$. One has

\begin{theorem}\label{thm nondeg fio}  \cite{Hor71,Hor85} Suppose that $C\subset  (T^*X\mzs)\times(T^*Y\mzs)$ is a canonical relation,
where $\dim(X)=n_1,\, \dim(Y)=n_2$, and $A\in I^{\meff-\frac{|n_1-n_2|}4}$ has a compactly supported Schwartz kernel. 
\smallskip

(i) If $C$ is nondegenerate, then $A:L^2_s(Y)\to L^2_{s-\meff}(X)$ for all $s\in\R$. 
Furthermore, the operator norm depends boundedly on a finite number of derivatives 
of the amplitude and phase function.
\smallskip 

(ii) If the spatial projections from $C$ to $X$ and to $Y$ are submersions
and, for some $l$,  the corank of $D\pi_L$ (and thus that of $D\pi_R$) is $\le l$ at all points of $C$,
then  $A:L^2_s(Y)\to L^2_{s-\meff-(l/2)}(X)$.
\end{theorem}
 
\medskip


\subsection{Frostman measures and $s$-energy}\label{subsec Frost}

Also   recall (see Mattila \cite{Mat95,Mat15})
that if $E\subset\R^d$ is a compact set  and $0<s<d$ satisfies $s<\hd(E)$, then
there exists a {Frostman measure} on $E$ relative to $s$: 
a probability measure $\mu$, supported on $E$,
satisfying the ball condition
\be\label{eqn ball}
\mu(B(x,\delta)\lesssim \delta^s,\, \forall\, x\in\R^d,\, 0<\delta<1,
\ee
 and of finite $s$-energy,

$$\int_E \int_E |x-y|^{-s}\, d\mu(x)\, d\mu(y) <\infty,$$
or equivalently, 
\be\label{eqn energy}
\int_E |\hat{\mu}(\xi)|^2\cdot |\xi|^{s-d}\, d\xi<\infty.
\ee
Since $\mu$ is of compact support, $\hat{\mu}\in C^\omega$ and thus  \eqref{eqn energy} 
implies that 
\be\label{eqn mu Sob}
\mu\in L^2_{(s-d)/2}(\R^d).
\ee
This also holds in the general setting of $E\subset X$,  
a compact subset of a $d$-dimensional manifold $X$ with $\hd(E)>s$.
\ms


\subsection{Tensor products of Sobolev spaces}\label{subsec Sob}

We need an elementary result on the tensor products of Sobolev spaces of negative order:

\begin{proposition}\label{prop Sob}
For $1\le j\le k$, let $X^j$ be a $C^\infty$ manifold of dimension $d_j$, and suppose that  $u_j\in L^2_{r_j,\, comp}\left(X^j\right),\, 1\le j\le k$,
with each $r_j\le 0$. Then the tensor product
$u_1\otimes \cdots \otimes u_k$ belongs to $ L^2_{r,\, comp}\left(X^1\times \cdots \times X^k\right)$,
for $r=\sum_{j=1}^k r_j$.
\end{proposition}
\noindent{\bf Proof.} Due to the compact support assumption, we can localize to a coordinate patch on each manifold,
reducing the problem to showing that 
$$L^2_{r_1}\left(\R^{d_1}\right)\otimes\cdots \otimes L^2_{r_k}\left(\R^{d_k}\right)\hookrightarrow 
L^2_r\left(\R^{\sum d_j}\right),$$
and this follows from the fact that each $\widehat{u_j}\left(\xi^j\right)\cdot \<\xi^j\>^{r_j}\in L^2\left(\R^{d_j}\right)$,
together with the lower bound $\Pi_{j=1}^k \<\xi^j\>^{r_j} \ge c\<\xi^1,\dots,\xi^k\>^r$ on $\R^{\sum d_j}$. Q.E.D.
\medskip


\subsection{Justification of density formula}\label{subsec density}

To justify \eqref{eqn pairing},  we argue as follows, restricting for simplicity the analysis to the case $k=3$ discussed in  Sec. \ref{sec three point},
when $\Phi:X^1\times X^2\times X^3\to\R^p$.
The proof  extends to $\Phi$ with codomain a general $T$ of dimension $p$ using local coordinates on $T$,
and also extends in a straightforward way to general $k$.

Without loss of generality, we consider $\sigma=(12|3)$.
For a $\chi\in C_0^\infty(\R^p)$  supported in a sufficiently small ball, $\chi\equiv 1$ near $\mathbf 0$,  and with $\int\chi\, d\vt=1$,
set $\chi_\epsilon(\vt):=\epsilon^{-p}\chi(\frac{\vt}{\epsilon})$ the associated approximation to the identity, 
which converges to $ \delta(\vt)$ weakly as $\epsilon\to 0^+$. 
Define $\rvte^{(12|3)}$ to be the operator with Schwartz kernel  
$$K_{\vt}^\epsilon(x^1,x^2;x^3):= \chi_\epsilon\left(\Phi\left(x^1,x^2,x^3\right)-\vt\right).$$
Then $\rvte^{(12|3)}( \mu_3)\in C^\infty(X^1\times X^2)$ and depends smoothly on $\vt$, and
thus we can represent the measure $\nu$ in \eqref{eqn pairing} as  the weak limit of absolutely continuous measures with smooth densities,
\be\label{eqn nu rep}
\nu(\vt)=\lim_{\epsilon\to 0^+} \nu^\epsilon(\vt):=\lim_{\epsilon\to 0^+}\langle\,\rvte^{(12|3)}\left(\mu_3\right),\mu_1\times \mu_2\rangle,
\ee
with $\nu$ having a density, which is in fact continuous in $\vt$, if the integral represented by the pairing converges.
Now, the operators $\rvte^{(12|3)}\in I^{-\infty}(C_{\vt}^{(12|3)})$, with symbols which converge in the Fr\'echet topology on the space of symbols  
as $\epsilon\to 0$ to the symbol of $\rvt$.
Since the singular limits $\rvt^{(12|3)}$ satisfy \eqref{eqn basic three} (for $\sigma=(12|3)$), so do the $\rvte^{(12|3)}$ uniformly in $\epsilon$. 
Hence, $\nu(\vt)$, being the uniform 
limit of smooth functions of $\vt$, is continuous. 
Furthermore, since $\epsilon^p\cdot \chi_\epsilon$ is bounded below by a constant times the 
characteristic function of the ball of radius $\epsilon$ in $\R^p$,    we have that
\be\label{eqn epsilon est}
\nu\left(B\left(\vt,\epsilon\right)\right):=
(\mu_1\times\mu_2\times\mu_3)\left(\left\{\left(x^1,x^2,x^3\right): \left|\Phi(x^1,x^2,x^3)-\vt\right|<\epsilon\right\}\right)\le C_\Phi\epsilon^p,
\ee
with constant $C_\Phi$ uniform in $\vt$, which was used in \eqref{eqn Leb} above.


\section{Proofs of theorems on 3-point configurations}\label{sec k is three}

 We are now able to prove the  theorems stated the Introduction that concern 3-point configurations:
 Thm. \ref{thm areas} about the areas  of  triangles  in $\R^2$ and  the radii of their circumscribing circles;
the three-dimensional case of Thm. \ref{thm volumes} regarding strongly pinned volumes of parallelepipeds;
and Thm. \ref{thm distsim} on ratios of pinned distances.
For all of these, we will show that Thm. \ref{thm threepoint} (ii) applies for appropriate choice of $\sigma$.
  \ms

\subsection{Areas  of triangles in $\R^2$.}\label{subsec areas proof}

We start with part (i) of Thm. \ref{thm areas}, on areas. The absolute value of the determinant is irrelevant for the 
conclusion of nonempty interior;
this will also be true for the other results where the configuration measurements have absolute values.
Additionally, the $1/2$ can be ignored. So, we start with the scalar-valued configuration function,
$$\Phi(x^1,x^2,x^3)
=\det\left[x^1-x^3,\, x^2-x^3\right]$$ 
on $\R^2\times\R^2\times \R^2$. 
Here, $d=2,\, p=1$, and we will show that,  using $\sigma=(12|3)$, that
the canonical relations $C^\sigma_t$ are (after  localizing) nondegenerate, so that 
Thm.  \ref{thm threepoint} (ii) yields a result for $\hd(E)>(2\cdot 2+1)/3=5/3$.
\ss
 
 We compute the gradient of $\Phi$
by noting that, on $\R^2\times\R^2$,
$$d_{u,v}\left(\det\left[u,v\right]\right)=\left(-v^\perp,u^\perp\right),$$
where $u^\perp=(-u_2,u_1)$ for $u=(u_1,u_2)$.
Hence, 
\be\label{eqn dPhi Athree}
d\Phi_{x^1,x^2,x^3}= \left(\left(x^3-x^2\right)^\perp,\, \left(x^1-x^3\right)^\perp,\, \left(x^2-x^1\right)^\perp\right).
\ee

Given a compact  $E\subset\R^2$ with $\hd(E)>5/3$, pick any $s$ with  $5/3<s<\nolinebreak\hd(E)$, 
and take $\mu$ to be a Frostman measure on $E$ with finite $s$-energy. 
Then we claim that one can find   points $x_0^1,\, x_0^2,\, x_0^3\in E$ and a $\delta>0$ such that
\footnote{Note that \eqref{eqn Athree generic2} just says that the $x_0^j$ belong to $supp(\mu)$ 
(which by Frostman's Lemma is $\subseteq E$ but can be a proper subset); however, for our purposes, it is useful to express this as \eqref{eqn Athree generic2}.}

\bea\label{eqn Athree generic1}
& &\det\left[x_0^1-x_0^3,\, x_0^2-x_0^3\right]\ne 0\hbox{ and }\\
& & \mu\left(B\left(x_0^j,\delta'\right)\right)>0,\, j=1,2,3, \forall \, 0<\delta'<\delta.\label{eqn Athree generic2}
\eea

To verify this, suppose not. Then, for every $x^1,\, x^2,\, x^3\in E$ and any $\delta>0$,  either 
\medskip

(i) $\Phi(x^1,x^2,x^3)=0$, i.e., $x^1,\, x^2,\, x^3$ are collinear, or 

(ii) for some $j=1,2$ or 3,  and some $\delta'<\delta$, $\mu(B(x^j,\delta'))=0$.
\medskip

Now, $Z_0=\{x\in\R^6:\, \Phi(x^1,x^2,x^3)=0\}$ is a five-dimensional algebraic variety.
Since $\mu\times\mu\times\mu$ has finite $3s$-energy, and  $3s>5$,  it follows that $(\mu\times\mu\times\mu)(Z_0)=0$,
and hence $(\mu\times\mu\times\mu)(E\times E\times E\setminus Z_0)=1$. 
We can in fact make this quantitative: Assuming without loss of generality  that $E$ is contained in the unit square centered at the origin,
 for $\epsilon>0$, let $\mathcal Z_\epsilon:=\left\{x\in\R^6: |x|<2 \hbox{ and }  \left|\Phi(x)\right|<\epsilon\right\}$.
 Then, since $Z_0$ is a rigid motion in $\R^6$ of the Cartesian product of $\R^2$ with a quadratic cone in $\R^4$,
 one sees that $\mathcal Z_\epsilon$ is covered by $\simeq \epsilon^{-5}$ balls of radius $\epsilon$ (away from the conical points),
 together with $\simeq \epsilon^{-1}$ balls of radius $\epsilon^{1/2}$ (covering a tubular neighborhood of the conical points). 
 Since $\mu$ satisfies the ball condition \eqref{eqn ball} on $\R^2$,  $\mu\times\mu\times\mu$ satisfies the corresponding condition on $\R^6$ with exponent $3s$
 and is thus dominated by $3s$-dimensional Hausdorff measure (up to a multiplicative constant). Thus,
 $$(\mu\times\mu\times\mu)(\mathcal Z_\epsilon)\lesssim \epsilon^{-5}\cdot \epsilon^{3s} + \epsilon^{-1}\cdot\epsilon^{3s/2}
 \lesssim \epsilon^{3s-5}\to 0 \hbox{ as } \epsilon\to 0.$$
Thus, if we define
$F_\epsilon:=E\times E\times E\setminus\mathcal Z_\epsilon$, 
which is compact, and ${\tilde\mu}_\epsilon:=(\mu\times \mu\times \mu)|_{F_\epsilon}$,
then ${\tilde\mu}_\epsilon(F_\epsilon)>1/2$ for $\epsilon$ sufficiently small.
By (ii) above, every $x\in F_\epsilon$ is in a $(\mu\times\mu\times\mu)$-null set which is also relatively open, 
the intersection of $F_\epsilon$ with a set 
of one of the three forms,
$$ B(x^1,\delta')\times \R^2\times \R^2, \quad \R^2 \times B(x^2,\delta')\times \R^2\hbox{ or } \R^2\times \R^2\times B(x^3,\delta').
$$
Since $F_\epsilon$ is compact, it is covered by a finite number of these, 
and hence it follows that $(\mu\times\mu\times\mu)(F_\epsilon)={\tilde\mu}_\epsilon(F_\epsilon)=0$. Contradiction. 
Hence, there exists an $x_0=(x^1_0,x^2_0,x^3_0) \in E\times E\times E$ such that \eqref{eqn Athree generic1},\eqref{eqn Athree generic2} hold. 
We now show that localizing near this base point allows us to apply Thm. \ref{thm threepoint}.
\medskip

Set $t_0=\Phi\left(x_0^1,\, x_0^2,\, x_0^3\right)\ne 0$; 
by continuity of $\Phi$  and relabelling there is a $\delta>0$ with $\Phi(x^1,x^2,x^3)\ne 0$ for   $x^j\in X^j:=B\left(x_0^j,\delta\right),\, j=1,2,3$. 
We claim that  for $\Phi|_{X^1\times X^2\times X^3}$ and $t$ close to $t_0$,  $(DF)_{(12|3)}$ is satisfied and and $C_t^{(12|3)}$ is nondegenerate, so that 
$\beta^{12,3}=0$ and  the last statement of Thm. \ref{thm threepoint} applies with $d=2$ and $p=1$; hence, if $\hd(E)>5/3$, then 
$\inter\left(A_3\left(E\right)\right)\ne\emptyset$.
\medskip

That $(DF)_{(12|3)}$ is satisfied is immediate, since all three components of $d_{x^1x^2,x^3}\Phi$ 
are nonzero on $X^1\times X^2\times X^3$ (the linear independence of the first two by \eqref{eqn Athree generic1} 
and the nonvanishing of third by that linear independence), which is an even stronger condition.
As for the canonical relations, one computes
$$C_t^{(12|3)}=\left\{\left(x^1,x^2,\theta\left(x^3-x^2\right)^\perp,\theta\left(x^1-x^3\right)^\perp;\, x^3,\theta\left(x^2-x^1\right)^\perp\right): 
(x^1,x^2,x^3)\in Z^{(12|3)}_t,\, \theta\ne 0\right\}.$$
Now, shrinking $\delta$ if necessary, for $\left(x^1,\, x^2,\, x^3\right)\in X^1\times X^2\times X^3$ and $t$ near $t_0$, 
for a smooth,  $X^1$-valued function $y^1(x^2,x^3,t)$ we can parametrize $Z_{t}^{(12|3)}$ by 
$$\left(y^1(x^2,x^3,t)+u(x^2-x^3),\, x^2,\, x^3\right),\quad (x^2,\, x^3)\in X^2\times X^3,\, u\in\R;$$
for example, one can take
$$y^1(x^2,x^3,t)=t \left|x^2-x^3\right|^{-2}\cdot (x^2-x^3)^\perp.$$
Thus, $(x^2,x^3,u,\theta)$ form coordinates on $C_t^{12,3}$, with respect to which
$$\pi_R\left(x^2,x^3,u,\theta\right)=\left(x^3,\,  \theta\left(x^2-x^1\left(x^2,x^3,u\right)\right)^\perp\right),$$
from which we see that
$$D_{x^3,u,\theta}\pi_R
=\left[
\begin{matrix}
\, I & 0 & 0 \cr
\, * & -\theta(x^2-x^3)^\perp & (x^2-x^1)^\perp
\end{matrix}
\right],
$$
which is of maximal rank since the last two columns are linearly independent.
Thus, $\pi_R$ is a submersion; by the general properties of canonical relations from 
Sec. \ref{subsec FIO}, $\pi_L$ is an immersion
and $C_t^{(12|3)}$ is nondegenerate. Q.E.D.
\bs

\subsection{Circumradii of triangles in $\R^2$.}\label{subsec circum proof}

We now turn to the proof of Thm. \ref{thm areas} (ii).
Changing the notation to denote the vertices of the triangle as $x,y,z\in\R^2$, 
 the circumradius $R(x,y,z)$ of $\triangle xyz$ is the distance from $x$ to the
intersection point of the perpendicular bisectors of $\overline{xy}$ and $\overline{xz}$.
For computational purposes, we  work with
\be\label{eqn twoRsquared}
\Phi(x,y,z):=2R^2(x,y,z)=\frac12\left(|y-x|^2+|z-x|^2+\frac{|z-x|^2\left(\left(y-x\right)\cdot\left(z-x\right)\right)^2}{\left(\left(y-x\right)^\perp\cdot
\left(z-x\right)\right)^2}\right),
\ee
with  the homeomorphism $r\to 2r^2$ of $\R_+$ of course preserving nonempty interior. 
With the same hypersurface $Z_0$ as in the proof of part (i) corresponds to degenerate triangles, 
$5/3<s<\hd(E)$ and Frostman measure $\mu$, as in (i)
we can find $x^0,\, y^0,\, z^0\in supp(\mu)$ such that 
all of the  components of $d\Phi$ are nonzero.
Picking $\sigma=(13|2)$, the incidence relation  $Z^\sigma_t=\left\{(x,z;y)\, :\, \Phi(x,y,z)=t\right\}$ satisfies  the 
double fibration condition $(DF)_\sigma$ 
for $t$ near $t_0=\Phi(x^0,y^0,z^0)$ and $x\in X^1,\, y\in X^2,\, z\in X^3$, neighborhoods of  $x^0,y^0,z^0$, resp. 
We can thus assume that $x_2\in\R,\, y\in\R^2,\, z\in\R^2$ form coordinates on $Z_t^\sigma$,
with $x_1$ a function of $x_2,y,z$, and rotating in $x$ if necessary, 
can further assume that $d_zx_1=0$ at $x^0,y^0,z^0$ and is therefore small nearby.
On $C_t^\sigma=(N^*Z_t^\sigma)'$, these together with the radial phase variable $\theta\in\R\setminus0$ are coordinates.
\ss

To show that the canonical relation $C_t^\sigma$ is nondegenerate, it suffices to show that $\pi_R: C_t^\sigma \to T^*X^2$ is a submersion.
Since the $y$ coordinate of $\pi_R(x_2,y,z,\theta)=y$, it suffices to show that 
$$rank\left[\frac{D\eta}{D(x_2,z,\theta)}\right]=2.$$
Setting
$$a=|y-x|^2,\quad b=(y-x)\cdot (z-x),\quad c=|z-x|^2,\quad d=(y-x)^\perp (z-y),$$
one calculates
\be\label{eqn Deta}
\frac{D\eta}{Dz}=\frac{cb}{d^2}\left(I-\frac{b}{d}J\right)-\left(1+\frac{cb}{d^2} -\frac{cb^2}{d^3}\right)
\left[\begin{matrix} 0 \\ d_zx_1\end{matrix}\right],
\ee
where $J=\left[\begin{matrix} 0 & -1 \\ 1 & 0 \end{matrix}\right]$ is the standard $2\times 2$ symplectic matrix, representing the $\perp$ map. 
The operator pencil $I-\lambda J\,  \lambda\in\R$, is nonsingular,
while the second term in \eqref{eqn Deta}  is small near $x^0,y^0,z^0$, and thus $D\eta/Dz$ is nonsingular, and $C_t^\sigma$ is nondegenerate.
Thus, as for areas, Thm.  \ref{thm threepoint} (ii) applies for $s>5/3$. Q.E.D.
\medskip

\subsection{Volumes of strongly pinned parallelepipeds in $\R^3$.} \label{subsec volumes 3d}

For  the proof of the $d=3$ case of  Thm. \ref{thm volumes}, 
the configuration function  $\Phi$ on  $\R^3\times\R^3\times\R^3$ is 
$$\Phi\left(x^1,x^2,x^3\right)=\det\left[x^1,\, x^2,\, x^3\right]=x^1\cdot\left(x^2\times x^3\right)= -x^2\cdot(x^1\times x^3)=x^3\cdot (x^1\times x^2).$$
We will show that  Thm. \ref{thm threepoint} (ii)  applies for $\sigma=(12|3)$, with $k=3,\, p=1$ and $d=3$, 
giving a positive result for $\hd(E)>(2\cdot 3 +1)/3=7/3$.
One computes
$$d\Phi_{x^1,x^2,x^3}=\left(x^2\times x^3, - x^1\times x^3, x^1\times x^2\right).$$
As in the previous proofs, given a compact  $E\subset\R^3$,  contained in the unit cube and  with $\hd(E)>7/3$,
pick $s$ with $7/3<s<\hd(E)$ and let   $\mu$ be a Frostman measure of finite $s$-energy. 
We claim  there exist  $x_0^1,\, x_0^2,\, x_0^3\in E$ and $\delta>0$ such that 
\medskip
\bea\label{eqn Vthree generic1}
& &x_0^1\times x_0^2\ne0,\, ,\, x_0^1\times x_0^3\ne 0,\, x_0^2\times x_0^3\ne 0,\hbox{ and }\\
& & \mu\left(B\left(x_0^j,\delta'\right)\right)>0,\, j=1,2,3, \forall \, 0<\delta'<\delta.\label{eqn Vthree generic2}
\eea

As before, we proceed with a proof by contradiction: suppose not. Then for every $x=(x^1,x^2,x^3)\in\R^9$ and $\delta>0$, either
\medskip

(i) at least one of $x^i\times x^j=0$, for some $1\le i<j\le 3$, or

(ii) for some $j=1,2$ or 3,  and some $\delta'<\delta$, $\mu(B(x^j,\delta'))=0$.
\medskip

On $\R^9$, $\mu\times\mu\times\mu$ has finite $3s$-energy and satisfies the ball condition with exponent $3s>7$.
For each $1\le i<j\le 3$ and $\epsilon>0$, 
$\mathcal W^{ij}_\epsilon:=\{x\in\R^9: |x^i\times x^j|<\epsilon\}$ is a tubular neighborhood of $\{x: x^i\times x^j=0\}$, 
a codimension two quadratic variety in $\R^9$ which 
is a rigid motion in $\R^9$ of the Cartesian product of $\R^3$ with a $4$-dimensional cone in $\R^6$. 
Following the analysis in the previous proof,
each of the $\mathcal W^{ij}_\epsilon$ can be covered by $\simeq \epsilon^{-7}$ balls of radius $\epsilon$ and $\epsilon^{-3}$ balls of radius $\epsilon^{1/2}$,
and thus 
$$(\mu\times\mu\times\mu)(\mathcal W^{ij}_\epsilon)\lesssim \epsilon^{3s-7} +\epsilon^{3s/2 -3}\lesssim \epsilon^{3s-7}\to 0\hbox{ as }\epsilon\to 0.$$
Hence, if we let $F_\epsilon=E\times E\times E\setminus \left(\cup_{i,j} \mathcal W^{ij}_\epsilon\right)$ 
and $\tilde\mu_\epsilon=(\mu\times\mu\times\mu)|_{F\epsilon}$,
then $F_\epsilon$ is compact and $\tilde\mu_\epsilon (F_\epsilon)>1/2$ for $\epsilon$ sufficiently small.
On the other hand, $F_\epsilon$ is covered by $\tilde\mu_\epsilon$-null and relatively open sets 
which are intersections of $F_\epsilon$ with sets 
of the three forms
$$ B(x^1,\delta')\times \R^3\times \R^3, \quad \R^3 \times B(x^2,\delta')\times \R^3\hbox{ or } \R^3\times \R^3\times B(x^3,\delta'),
$$
and the compactness of $F_\epsilon$ leads to a contradiction. Hence, we can find $x_0^1,x_0^2,x_0^3$ and $\delta$ such that \eqref{eqn Vthree generic1},
\eqref{eqn Vthree generic2} hold. Further restricting $\delta$ if necessary, we can assume that
$x^1\times x^2,\, x^1\times x^3$ and $x^2\times x^3$ are $\ne 0$.  
for all $x^j\in B\left(x_0^j,\delta\right)=:X^j,\, j=1,2,3$. 

Restricting $\Phi$ to $X^1\times X^2\times X^3$, $(DF)_{(12|3)}$ is satisfied; in fact all three components of $d\Phi$ are nonzero.
On the incidence relation $Z_t^{(12|3)}$, we can take as coordinates $x^2,\, x^3$ 
and $\vec{u}=(u_2,u_3)\in\R^2$, solving for $x^1$ with
$$x^1=y^1(x^2,x^3,t) + u_2x_2+u_3x_3,$$
for some smooth function $y^1$. Thus,

\bast
C_t^{(12|3)}&=&\big\{\left(x^1,x^2,\theta\left(x^2\times x^3\right), -\theta\left(x^1\times x^3\right);\, x^3, -\theta\left(x^1\times x^2\right)\right) \\
& & \qquad\qquad\qquad : 
(x^2,x^3)\in X^2\times X^3,\, \vec{u}\in\R^2,\, \theta\ne 0\big\}.
\east
Then $\pi_R$ is a submersion since, with $\xi^3=-\theta\left(x^1\times x^2\right)$, $D_{x^2,u_3,\theta}\xi^3$  is surjective:  one has 
$D\xi^3(\partial_\theta)=x^1\times x^2,\, D\xi^3(\partial_{u_3})=\theta\left(x^2\times x^3\right)$ and the range of $D_{x^2}\xi^3$ is $\left(x^1\right)^\perp$;
together, these span all of the $\partial_{\xi^3}$ directions.

Since $\pi_R$ is a submersion, $C_t^{(12|3)}$ is nondegenerate, and  
Thm. \ref{thm threepoint} (ii) applies, this time with $d=3$ and $p=1$; hence, if $\hd(E)>7/3$, then 
$\inter\left(V_3^0\left(E\right)\right)\ne\emptyset$. Q.E.D.

\begin{remark}\label{rem two d pinned}
 The proof for  $d\ge 4$ will be presented in Sec. \ref{subsec pinned proof} below. 
On the other hand, Thm. \ref{thm volumes} does not give a positive result for   pinned volumes (areas)  in two dimensions,
$$V_2^0(E):=\left\{\det\left[x^1,\, x^2\right]:\, x^1,x^2\in E\right\}.$$
In fact, since this concerns a 2-point configuration, it would already fall under the framework of
\cite{GIT19}; however, the projections $\pi_L,\, \pi_R$ from the canonical relation to $T^*\R^2$ both drop rank by 1 everywhere,  
resulting in a loss of $\beta^{(1|2)}=1/2$ derivatives. Hence, $\mathcal R_t^{(1|2)}\in I^{-\frac12}\left(C_t^{(1|2)}\right)$ is not smoothing on $L^2$-based Sobolev spaces,
and the FIO approach to configuration problems does not imply a result in this case.
\end{remark} 

\subsection{Ratios of pinned distances}\label{subsec distsim}

For $E\subset\R^d,\, d\ge 2$, we prove Thm. \ref{thm distsim} concerning the  set defined in \eqref{eqn distsim}. 
On $\left(\R^d\right)^3$, let
$$\Phi(x^1,x^2,x^3)=\frac{|x^1-x^3|}{|x^1-x^2|}.$$ 
We show that, after suitable localization,  Thm. \ref{thm threepoint} (ii), with $k=3,\, p=1$,  applies for $\sigma=(12|3)$, 
implying a nonempty interior result  when $\hd(E)>(2d +1)/3$.

One computes,
$$d\Phi(x^1,x^2,x^3)=\left| x^1-x^2\right|^{-2} \left( (x^2-x^3), -(x^1-x^2),-(x^1-x^3)\right).$$
Let $\hd(E)>(2d+1)/3$ and $\mu$ be a Frostman measure on $E$ of finite $s$-energy for some $(2d+1)/3 <s <\hd(E)$.
Then $\mu\times\mu\times\mu$ is dominated by $3s$-dimensional Hausdorff measure, and $3s>2d+1>2d$.
Since $\{(x^1-x^2)(x^1-x^3)(x^2-x^3)=0\}$ is a union of three $2d$ dimensional planes, as above one can show  
that there exist $x_0^1,x_0^2,x_0^3\in supp(\mu)$ such that $x_0^i-x_0^j\ne 0,\, i\ne j$. 
Taking $X^j=B(x_0^j,\delta)$ for suitably small $\delta$, all three components of $d\Phi$ are nonzero and,
setting $t_0=\Phi(x_0^1,x_0^2,x_0^3)$, 
the double fibration condition $(DF)_{\sigma}$ is satisfied by $Z_t^\sigma$ for $t$ close to $t_0$.

On $Z_t^\sigma$, solving for $x^3=x^1-t|x^1-x^2|\omega$ we can take as coordinates $(x^1,x^2,\omega)\in \R^d\times\R^2\times\sd$
. Then, the canonical relation of $\zts$ is
\bast
C_t^\sigma&=&(N^*\zts)' \\
&=&\left\{ \left(\cdot,\cdot,\cdot,\cdot\, ;\,  x^1-t|x^1-x^2|\omega, \theta t|x^1-x^2|\omega\right)\, 
:\, x^1,x^2\in\R^d,\, \omega\in\sd,\,\theta\ne 0 \right\},
\east
where we have suppressed the $T^*\left(X^1\times X^2\right)$ components as irrelevant for 
analyzing $\pi_R$. 
One easily sees that the projection from $C_t^\sigma$ to $T^*X^3$ is a submersion, so that Thm. \ref{thm threepoint} (ii) applies as claimed.

 
 \section{$k$-point configuration sets, general $k$}\label{sec bigger than three}
 
 To describe results for general $k$-point configurations, let $X^i,\, 1\le i\le k,$ and $T$, be  
 smooth manifolds of dimensions $d_i$ and $p$, 
 resp. We sometimes denote \linebreak $X^1\times\cdots\times \nolinebreak X^k$ by $X$, and set $\dt:=\dim(X)=\sum_{i=1}^k d_i$.
 
 \begin{definition}\label{def k Phi config} 
  Let $\Phi\in C^\infty(X,T)$. 
Suppose that   $E_i\subset X^i,\, 1\le i\le k$, are compact sets.
Then the {\it $k$-configuration set of the $E_i$} defined by $\Phi$ is
\be\label{def config general}
\Delta_\Phi\left(E_1,E_2,\dots,E_k\right):=\left\{\Phi\left(x^1,\dots,x^k\right): x^i\in E_i,\, 1\le i\le k\right\}\subset T.
\ee
\end{definition}

\medskip

We want to find sufficient conditions on  $\hd(E_i)$ ensuring that  $\Delta_\Phi\left(E_1,E_2,\dots,E_k\right)$
has nonempty interior.
To this end, now suppose that $\Phi:X\to T$ is an submersion, so that for each $\vt\in T$, $\zvt:=\Phi^{-1}(\vt)$ is a smooth, 
codimension $p$ submanifold of $X$, 
and these vary smoothly with $\vt$. 
For each $\vt$, the measure
 \be\label{def delta k}
\lambda_\vt:= \delta\(\Phi\(x^1,\dots,x^k\)-\vt\)
 \ee
 is a smooth density on $\zvt$; using local coordinates on $T$, one sees that this can be represented as an oscillatory integral
of the form
 $$\int_{\R^p} e^{i\[\sum_{l=1}^p\(\Phi_l(x^1,\dots,\, x^k)-\vt_l\)\theta_l\]} a(\vt)1(\theta)\, d\theta,$$
 where the $a(\cdot)$ belongs to a partition of unity on $T$. 
 Thus, $\lambda_\vt$ is a Fourier integral distribution on $X$; 
 in H\"ormander's notation \cite{Hor71,Hor85},
 \be\label{eqn fid}
 \lambda_\vt\in I^{(2p-\dt)/4}(X;N^*\zvt),
 \ee
  where $N^*\zvt\subset T^*X\setminus 0$ is the conormal bundle of $\zvt$ and
  the value of the  order  follows
from the amplitude having order zero and the numbers of phase variables and spatial variables  being $p$ and $\dt$, resp., 
so that the order is $m:=0+p/2-\dt/4$.
\ms

As   in the analysis of $3$-point configurations in Sec. \ref{sec three point}, 
we separate the variables $x^1,\,\dots, x^k$ into groups on the left and right,
  associating to $\Phi$ a collection of families of generalized Radon transforms  
indexed by the nontrivial partitions of $\{1,\dots,k\}$, 
with each family then 
depending on  the parameter $\vt\in T$.
Write such a partition   as $\sigma=(\sigma_L\, |\, \sigma_R)$,
with $|\sigma_L|,\, |\sigma_R|>0,\, |\sigma_L|+|\sigma_R|=k$, 
and let $\pk$ denote the set of all such partitions.
We will use $i$ and $j$ to refer to elements of $\sigma_L$ and $\sigma_R$, resp.
Define $\dl=\sum_{i\in \sigma_L} d^i$ and $\dr=\sum_{i\in \sigma_R} d^i$, so that $\dl+\dr=\dt$.
\ms

For each $\sigma\in\pk$, $\sigma_L=\left\{i_1,\dots,i_{|\sigma_L|}\right\}$ and $\sigma_R=\left\{j_1,\dots, j_{|{\sigma_R}|}\right\}$,
where we may assume that $i_1<\cdots i_{|\sigma_L|}$ and $j_1<\cdots< j_{|\sigma_R|}$.
with a slight abuse of notation
we still refer to as $x$ the permuted version as $x$,
$$x=\left(x_L;x_{R}\right):=\(x^{i_1},\dots,x^{i_{|\sigma_L|}};\,  x^{j_1},\dots,x^{j_{|\sigma_R|}}\).$$

Write the  corresponding reordered  Cartesian product as
$$X_L \times X_{R} := \(X^{i_1} \times \cdots \times X^{i_{|\sigma_L|}}\) \times
\(X^{j_1} \times \cdots \times X^{j_{|\sigma_R|}}\);$$
again by abuse of notation, we sometimes still refer to this as $X$.
The dimensions of the two factors are $dim(X_L)=\dl$ and $dim(X_R)=\dr$, resp.
The choice of $\sigma$  also  defines a permuted version of each $\zvt$, 
\be\label{def zsigma}
\zvts:=\left\{ \(x_L; x_{R}\):\, \Phi\(x\)=\vt  \right\}\subset X_L \times X_{R},
\ee
with spatial projections to the left and right,  $\pi_{X_L}:\zvts\to X_L$ and $\pi_{X_R}:\zvts\to X_{R}$.
The  integral geometric double fibration condition extending \eqref{def DF} to general $k$ is
the requirement, 
\be\label{def DF general}
 (DF)_{\sigma}\qquad \pi_{L}:\zvts\to X_L\hbox{ and }\pi_{R}:\zvts\to X_{R}\hbox{ are submersions.}
\ee 
(Note that  $(DF)_\sigma$ can only hold for a given $\sigma$ if $p\le \min\left(\dl, \dr\right)$.)
\medskip

If $(DF)_\sigma$  holds, then the generalized Radon transform $\rvts$,
defined weakly by
$$\rvts f(x_L)=\int_{\left\{x_{R}:\, \Phi\left(x_L,x_{R}\right)=\vt\right\}} f(x_{R}),$$
where the integral is with respect to the surface measure induced by $\lambda_\vt$ on the codimension $p$ submanifold
$\left\{x_{R}:\, \Phi\left(x_L,x_{R}\right)=\vt\right\}=\left\{x_{R}: \left(x_L,x_{R}\right)\in \zvts\right\} \subset X_{R}$,
extends from mapping $\mathcal D\left(X_{R}\right)\to \mathcal E\left(X_L\right)$ to
$$\rvts:\mathcal E'(X_{R})\to \mathcal D'(X_L).$$
Furthermore,  
\be\label{def ctsigma}
C_{\vt}^\sigma:=\(N^*\zvts\)'=\left\{\(x_L,\xi_L;\, x_{R},\xi_{R}\):
\, \(x_L,x_{R}\)\in\zvts,\, \(\xi_L,-\xi_{R}\)\perp T\zvts \right\}
\ee
is contained in $\(T^*X_L\setminus 0\)\times \(T^*X_R\setminus 0\)$. 
Thus, $\rvts$ is an FIO,  $\rvts\in I^m\(X_L, X_R; C_{\vt}^\sigma\)$,
where the order $m$ is determined by H\"ormander's formula, $m=0+p/2-\dt/4$. 
Given the possible difference in the dimensions of $X_L$ and 
$X_R$,  due to the clean intersection calculus it is useful to express $m$  as
$$m=\meff^\sigma-\frac14\left|\dl-\dr\right|,$$
where the {\it effective} order of $\rvts$  is defined to be
\be\label{eqn msigma}
\meff^\sigma:=\(2p-\dt+\left|d^\sigma-d^{\hs}\right|\)/4=\left(p-\min\left(\dl,\dr\right)\right)/2.
\ee
As recalled in Sec. \ref{subsec FIO}, if $C_{\vt}^\sigma$  is a nondegenerate canonical relation, i.e., the  cotangent space projections
$\pi_L: C_{\vt}^\sigma \to T^*X_L$
and 
$\pi_R:C_{\vt}^\sigma \to T^*X_R$
have differentials of  maximal rank\footnote{This is a structurally stable condition, so that if  $C_{\vt_0}^\sigma$ is nondegenerate, 
then $C_{\vt}^\sigma$ is nondegenerate for all $\vt$ 
in some neighborhood of $\vt_0$.}, then 
$$\rvts:L^2_r\left(X_R\right)\to L^2_{r-{\meff}^\sigma}\left(X_L\right).$$
As in the result concerning 3-point configurations  in Thm. \ref{thm threepoint}, 
it is natural to express the estimates for possibly degenerate FIO  in terms of possible losses relative to the optimal estimates.
As for $k=3$, our basic  assumption is  that,  for at least one $\sigma$, the double fibration condition \eqref{def DF general} is satisfied and 
there  is a known $\beta^\sigma\ge 0$ such that, for all $r\in\R$, 
\be\label{eqn basic k}
\rvts: L^2_r\left(X_R\right) \to L^2_{r-\meff^\sigma-\beta^\sigma}\left(X_L\right)\hbox{ uniformly for }\vt\in T.
\ee

Now  suppose that,  for $1\le i\le k$,   $E_i\subset X^i$ are compact sets.
Our goal is  to find  conditions on the $\hd(E_i)$ ensuring that 
$\Delta_\Phi\left(E_1,E_2,\dots,E_k\right)$
has nonempty interior in $T$.
For each $i$, fix an $s_i<\hd(E_i)$ 
and  a Frostman measure  $\mu_i$ on $E_i$ of finite $s_i$-energy.
Define measures  
$$\mu_L := \mu_{i_1}\times \cdots \times \mu_{i_{|\sigma|}}\hbox{ on }X_L\hbox{ and }
\mu_R := \mu_{j_1}\times \cdots \times \mu_{j_{|\hs|}}\hbox{ on }X_R.$$
By \eqref{eqn mu Sob}, each $\mu_i\in L^2_{(s_i-d_i)/2}(X^i)$,
so that Prop. \ref{prop Sob} implies that 
  $\mu_L\in L^2_{r_L}\left(X_L\right)$ and
 $\mu_R\in L^2_{r_R}\left(X_R\right)$, 
 where $r_L=\frac12\sum_{l=1}^{|\sigma_L|} \left(s_{i_l}-d_{i_l}\right)$ and
 $r_R=\frac12\sum_{l=1}^{|\sigma_R|} \left(s_{j_l}-d_{j_l}\right)$, resp.
The analogue of the representation formula \eqref{eqn nu rep} for $\nu(\vt)$,
justified by a minor modification of the $k=3$ case in Sec. \ref{subsec density},
 is
 \be\label{eqn pairing k}
 \nu(\vt)=\<\rvts\left(\mu_R\right),\mu_L\>.
 \ee

Our basic assumption, that \eqref{eqn basic k} holds for the $\sigma$ in question, then implies 
that $\rvts\left(\mu_R\right)\in L^2_{r_R-\meff^\sigma-\beta^\sigma}\left(X_L\right)$.
 Since $\mu_L\in L^2_{r^\sigma}\left(X_L\right)$, the pairing in  \eqref{eqn pairing k} is bounded,
 and yields a continuous function of $\vt$, if
 \be\label{eqn first ineq}
 r_R-\meff^\sigma-\beta^\sigma + r_L \ge 0.
 \ee
 Noting that 
 $$r_L+r_R=\frac12\left[\left(\sum_{i=1}^k s_i\right) - d^{tot}\right],$$
 and using \eqref{eqn msigma}, we see that \eqref{eqn first ineq} holds iff
 \bast
 \sum_{i=1}^k s_i &\ge & d^{tot} +2\left(\meff^\sigma + \beta^\sigma\right) \\
 & = & d^{tot} + p - \min\left(d_L,d_R\right) + 2\beta^\sigma \\
 & = & \max\left(d_L,d_R\right) + p + 2\beta^\sigma.
  \east
  
  Optimizing over all $\sigma \in \pk$, we obtain the analogue of Thm. \ref{thm threepoint} for $k$-point configuration sets:

\begin{theorem} \label{thm kpoint}
(i) With the notation and assumptions as above, define
$$s_\Phi=\min \left (\max\left(d_L,d_R\right) + p + 2\beta^\sigma \right),$$
where the $\min$ is taken over those  $\sigma\in\pk$ for which both  the double fibration condition  
\eqref{def DF general} 
and the uniform boundedness of the generalized Radon transforms $\rvt^\sigma$ with loss of $\,\le \beta^\sigma$ derivatives 
\eqref{eqn basic k} hold.
\ms

Then, if $E_i\subset X^i$, $1\le i\le k$, are compact sets with $\sum_{i=1}^k\hd(E_i)>s_\Phi$, 
it follows that $\inter\left(\Delta_\Phi\left(E_1,E_2, \dots , E_k\right)\right)\ne\emptyset$.
\bs

(ii) In particular, if $X^1=\cdots=X^k=X_0$, with $dim(X_0)=d$,
and if $E\subset X_0$ is compact with 
\be\label{eqn equi k}
\hd(E)>\frac{1}{k}\left[ \min \left(\max\left(d_L,d_R\right) + p\right) \right],
\ee
where the minimum is taken over all
$\sigma\in\pk$ such that the canonical relations $C_{\vt}^\sigma$ are 
nondegenerate,
then
 $\inter\left(\Delta_\Phi\left(E,E,\dots, E\right)\right)\ne\emptyset$.  
\end{theorem}


\section{Proofs of theorems for $k$-point configurations, $k\ge 4$}\label{sec k is four or more}

We now use Thm. \ref{thm kpoint} to prove Thm. \ref{thm cross} on cross ratios of four-tuples of points in $\R$,
Thm.  \ref{thm volumes} concerning strongly pinned volumes of parallelepipeds generated by $d$-tuples of points in $\R^d$ for $d\ge 4$,
Thm. \ref{thm two areas} on pairs of areas of triangles in $\R^2$, 
Thm. \ref{thm prod of diff} on dot products of differences
and 
Thm. \ref{thm bilinear forms} on generalized  sum-product sets.

\subsection{Cross ratios on $\R$}\label{subsec cross}

We prove Thm. \ref{thm cross} on the set of cross ratios of four-tuples of points in a set $E\subset\R$.
This will be an application of the second part of Thm. \ref{thm kpoint}, with $d_j=1,\, 1\le j\le 4$. 
Since each of the variables in one dimensional,
we use subscripts rather than superscripts.
Thus, set $\Phi:\R^4\to\R$, 
$$\Phi(x_1,x_2,x_3,x_4) = \left[x_1,x_2;x_3,x_4\right]
=\frac{(x_1-x_3)(x_2-x_4)}{(x_1-x_4)(x_2-x_3)},$$
and introduce the notation 
$$x_{ij}=x_i-x_j,\quad 1\le i<j\le 4.$$
Then one computes
$$d\Phi(x_1,x_2,x_3,x_4)=\left(x_{14}x_{23}\right)^{-2}\big(x_{23}x_{24}x_{34},\, -x_{13}x_{14}x_{34},\, x_{14}x_{24}x_{34},\, -x_{12}x_{13}x_{23}\big).$$
Given a compact $E\subset\R$ with $\hd(E)>3/4$, let $s$ be such that $3/4<s<\hd(E)$ and $\mu$ be a Frostman measure on $E$ of finite $s$-energy.
We claim that 
\bea\label{eqn cross generic}
\left(\exists\, x_1^0,\, x_2^0,\, x_3^0,\, x_4^0\in supp\left(\mu\right)\right)\hbox{ s.t. }
x_i^0-x_j^0\ne 0,\quad\hbox{ for all }1\le i<j\le 4,
\eea
so that all four components of $d\Phi$ are nonzero at $x^0:=(x_1^0,x_2^0,x_3^0,x_4^0)$.
Arguing as in the proofs for $3$-point configurations in Sec. \ref{sec k is three}, 
set 
$$\mathcal Z=\left\{x\in\R^4\, :\, \prod_{1\le i<j\le 4} x_{ij}=0\right\},$$
on the complement of which all of the components of $d\Phi$ are nonzero.
Noting that $\mathcal Z$ is a union of hyperplanes,  $\hd(\mathcal Z)=3$; 
thus, since $\mu\times\mu\times\mu\times\mu$ is dominated by $4s$-dimensional Hausdorff measure and $4s>3$,
  $\left(\mu\times\mu\times\mu\times\mu\right)(\mathcal Z)=\nolinebreak0$.
By a slight variant of the reasoning in Thms.  \ref{thm areas} and Thm. \ref{thm volumes}, one obtains \eqref{eqn cross generic}.
This actually shows that the conditions in  \eqref{eqn cross generic} hold for a set of full $\mu\times\mu\times\mu\times\mu$ measure, which we use below.
\medskip

Now let $\sigma=(12|34)$.
Setting $t^0=\Phi(x^0)$, and taking the $X^j$ to be sufficiently small neighborhoods of $x_j^0,\, 1\le j\le 4$, 
it follows that for $t$ close to $t^0$ the double fibration condition \eqref{def DF general} holds on 
$$\zts:=\left\{x\,:\, \Phi(x)=t\right\}\subset X:=(X^1\times X^2)\times (X^3\times X^4) =:X_L\times X_R.$$
Since $d_{x_1}\Phi\ne 0$, on $Z_t^\sigma$ we can solve for $x_1$ as a function of $x_2,x_3,x_4$ and the parameter $t$
(possibly again reducing the size of the neighborhoods of the $x_j^0$),
$$x_1=y^1(x_2,x_3,x_4,t).$$

Furthermore, using the fact that the factor $\left(x_{14}x_{23}\right)^{-2}$ in all of the terms of $d\Phi$ can be absorbed into the radial scaling factor,
as in \eqref{def ctsigma} we define
\bast
C_{t}^\sigma:=\(N^*\zts\)'&=&\big\{\left(*,\, *; x_3,\, x_4, \theta x_{14}x_{24}x_{34},\, -\theta x_{12}x_{13}x_{23}\right)\, \\
& &\qquad  :\,( x_2,x_3,x_4)\in X^2\times X^3\times X^4,\, \theta\ne 0\big\},
\east
where we have suppressed the $T^*(X^1\times X^2)$ components on the left.
We claim that this is a local canonical graph, so that the family of FIOs lose no derivatives ($\beta^\sigma=0$); 
the desired result then follows from Thm. \ref{thm kpoint} (ii),
using  $d_L=d_R=2,\, k=4$ and $p=1$ in \eqref{eqn equi k}. 
 
 \medskip
 
 To see that $C_t^\sigma$ is a local canonical graph, it suffices to show that the differential of $\pi_R:C_t^\sigma\to T^*X_L$ 
 has rank 4. Due to the $(x_3,x_4)$ in the spatial variables, this is equivalent to showing that 
 $$\frac{D(\xi_3,\xi_4)}{D(\theta,x_2)}=\left[
 \begin{matrix}x_{14}x_{24}x_{34}, & \theta\left(x_{34}\left(x_{14}+x_{24}y^1_{x_2}\right)\right) \\
 -x_{12}x_{13}x_{23}, & -\theta\left(-x_{13}x_{23}+x_{12}x_{13}+x_{23}\left(x_{13}y^1_{x_2}+x_{12}y^1_{x_2}\right)\right)
 \end{matrix}
 \right]
 $$
 is nonsingular. However, the determinant of this is an algebraic function, not identically vanishing on $Z_t^\sigma$, so its zero variety is three-dimensional and thus a null set with respect to $4s$-Hausdorff measure, and hence with respect to $\mu\times\mu\times\mu\times \mu$.
 Thus,  choosing our basepoint $x^0$, and then shrinking the $X^j$ suitably,  to avoid this, ensures that $C_t^\sigma$ is a local canonical graph, finishing the proof of Thm. \ref{thm cross}.
 

\subsection{Strongly pinned volumes in $\R^d,\, d\ge 4$}\label{subsec pinned proof}

 With Thm. \ref{thm kpoint} in hand, we now prove Thm.  \ref{thm volumes} concerning pinned volumes in $\R^d$ for $d\ge 4$, 
 following the lines of the proof  for $d=3$ in Sec. \ref{sec k is three}.  
 On $\left(\R^d\right)^d$, let $\Phi\left(x^1,\dots,x^d\right)=\det\left[x^1, x^2,\dots, x^d\right]$.
 We will show that for $\sigma=\left(12\dots \left(d-1\right)|d\right)$,  some $t_0\ne 0$ and with the domain of $\Phi$ suitably localized,  
 condition $(DF)_\sigma$ is satisfied and  the canonical relation $C_{t_0}^{\sigma}$ is nondegenerate;
 these conditions then hold for all $t$ near $t_0$ by structural stability of submersions.
Using just this  $\sigma$,  applying \eqref{eqn equi k} with
$d_L=d(d-1)>d_R=d$, 
 $p=1$ and $\beta_\Phi=0$ shows that if $\hd(E)>\left(1/d\right)\left(d\left(d-1\right)+1\right)=d-1+\left(1/d\right)$
then $\inter\left(\Delta_\Phi\left(E,\dots,E\right)\right)\ne\emptyset$, proving Thm. \ref{thm volumes}.
\smallskip
 
 To verify the claims for  $\sigma$, we start by noting that
 $$d\Phi=\left(\bx^{(1)}, -\bx^{(2)},\dots,(-1)^{d-1}\bx^{(d)}\right),$$
 where
 $$\bx^{(j)}:=*\left(x^1\wedge x^2\wedge \cdots \wedge x^{j-1}\wedge x^{j+1}\wedge\cdots \wedge x^d\),$$
 where $*$ is the Hodge star operator, which is an isomorphism $*:\Lambda^{d-1}\R^d\to \R^d$.
 As in the proof for $d=3$, note that
 if $d-1+\left(1/d\right)<s<\hd(E)$ and $\mu$ is a Frostman measure on $E$ of finite $s$-energy, 
one {can find} $x_0^1,\, x_0^2,\, \dots,\, x_0^d\in supp(\mu)$ and $\delta>0$ such that
$\bx^{(j)}\ne 0,\, 1\le j\le d$,
whenever $x^j\in B\left(x_0^j,\delta\right)=:X^j,\, 1\le j\le d$.
This follows by a straight-forward modification of the argument in Sec. \ref{subsec volumes 3d}.

For $1\le j\le d$, each variety $\mathcal W^{(j)}:=\left\{x\in \R^{d^2}: \bx^{(j)}=0\in\R^d\right\}$  
is codimension $(d-1)$, and thus their union is a null set with respect to 
$\otimes^d \mu$, since $sd>d^2-(d-1)$.
Restricting $\Phi$ to $X^1\times \cdots \times X^d$, $(DF)_{\sigma}$ is satisfied. In fact, each of the  components of $d\Phi$,
$$ d_{x^j}\Phi|_{x_0}= (-1)^j\bx^{(j)},
$$
 is nonzero, since when paired against $x_0^j$ it gives $\Phi(x_0^1,\dots,x_0^d)\ne 0$. 
 \smallskip

Thus, for $\sigma=\left(12\dots \left(d-1\right)|d\right)$ and $t$ close to $t_0$, as 
coordinates on the incidence relations $Z_t^{\sigma}$ we can take  
$(x^2,\, \dots,\, x^d)\in X^2\times\cdots \times X^d$ and $\vec{u}=(u_2,\, \dots,\, u_d)\in\R^{d-1}$, 
with $x^1$ determined by 
$$x^1=y^1\left(x^2,\, \dots,x^d,t\right) + u_2x^2+\cdots+u_dx^d,$$
for some smooth function $y^1$, since the perturbations of any specific $x^1$ 
that preserve $\det\left[x^1, x^2,\dots, x^d\right]$
are arbitrary translates in the directions spanned by $x^2,\dots,x^d$. 
Furthermore, by translating by a constant in the $s$ variables, 
we can assume that at the base point,
\be\label{eqn Dxtwo fact}
D_{x^2}y^1(\vec{x}_0,t_0)=0.
\ee
Thus,
in $T^*\left(X^1\times\cdots\times X^{d-1}\right)\times T^*X^d$, 
\bast
C_t^\sigma\!\!\!&=&\!\!\!\big\{\big(\,\cdot,\, \cdot\, ; x^d,\, \pm * \theta
\left[\left(y^1\left(x^2,\, \dots,x^d,t\right) \wedge x^2\wedge\cdots\wedge x^{d-1}\right)+ (-1)^d \left(\, u_d x^2\wedge \cdots \wedge x^d\right)\right] \big)\\
& &\qquad\qquad\qquad  : (x^2,\dots,x^d)\in X^2\times\cdots X^d,\, \vec{u}\in\R^{d-1},\, \theta\in \R\setminus 0  \big\},
\east
where the first  entries, giving the coordinates in $T^*\left(X^1\times\cdots X^{d-1}\right)$, have been suppressed because
they are not needed to study $\pi_R$. In the last, i.e., $\xi^d$, entry, we have used 
\bast
\big(y^1(x^2,\, \dots,x^d,t) &+& u_2x^2+\cdots+u_dx^d\big)\wedge\, x^2\wedge\cdots\wedge x^{d-1}\\
&= &y^1(x^2,\, \dots,x^d,t) \wedge\, x^2\wedge\cdots\wedge x^{d-1}\\
& & + (-1)^d \, u_d\, x^2\wedge \cdots \wedge x^d.
\east
We claim that $\pi_R: C_t^\sigma\to T^*X^d$ is a submersion, which, as described in Sec. \ref{subsec FIO},
 then implies that $C_t^\sigma$ is nondegenerate and thus $\sigma$ is one of the competitors in \eqref{eqn equi k}.
Note that
$$D_{x^d}\pi_R={\mathbf I_d}\oplus\left( \pm *\theta\left[ D_{x^d}y^1\wedge x^2\wedge\cdots\wedge x^{d-1} 
+\left(-1\right)^d u_d\, x^2\wedge\cdots \wedge x^{d-1}\wedge I_d\right]\right)$$
while, for $2\le j\le d-1$,
\be\label{eqn Dxtwo}
D_{x^j}\pi_R={\mathbf 0}\oplus \left( \pm *\theta\left[ D_{x^j}y^1\wedge x^2\wedge\cdots\wedge x^{d-1} 
+\left(-1\right)^d u_d\, x^2\wedge\cdots\wedge I_d \wedge\cdots\wedge x^{d}\right]\right).
\ee
 
 Due to the form of $D_{x^d}\pi_R$, it suffices to show that  $D_{x^2}\pi_R$ has rank equal to $d$.
Since $\theta\ne 0$ and $*$ is an isomorphism, we can ignore the $\pm *\theta$ and work directly in  the
 $d$-dimensional vector space $\Lambda^{d-1}\R^d$. At $x_0$, the expression in square brackets in \eqref{eqn Dxtwo}
 equals $(-1)^du_dI_d\wedge\left(x^3\wedge \cdots \wedge x^d\right)$ due to \eqref{eqn Dxtwo fact}.
 Since $x^3\wedge \cdots \wedge x^d\in\Lambda^{d-2}\R^d-\{0\}$, 
 this last map is an isomorphism $\R^d\to \Lambda^{d-1}\R^d$, thus has rank $d$, finishing the proof.

  
  \subsection{Pairs of areas in $\R^2$}\label{subsec pairs}

 We now prove Thm. \ref{thm two areas} concerning the set  of  pairs of areas of triangles generated by 4-tuples of  points  in a compact $E\subset\R^2$. 
Here, $d=2,\, k=4$ and $p=2$.
 On $\left(\R^2\right)^4$, let
 \bast
 \Phi\left(x^1,x^2,x^3,x^4\right)&=& \left(\det\left[x^1-x^4,x^2-x^4\right],\, \det\left[ x^2-x^4,x^3-x^4\right]\right)\\
 &=& \big((x^1-x^4)\cdot(x^2-x^4)^{\perp}, (x^2-x^4)\cdot(x^3-x^4)^{\perp}))
 \east
We will show that, for $\sigma=(13|24)$, although  $C_t^\sigma$ is degenerate, the projections $\pi_L,\, \pi_R$
 drop rank by at most $1$ everywhere,
 and therefore, by Thm. \ref{thm nondeg fio}(ii), there is a loss of at most $\beta^\sigma=1/2$ derivative. 
 Here, $d_L=d_R=4$, so Thm. \ref{thm kpoint} implies that for $$4\, \hd(E)> \max(d_L,d_R)+p+2\beta_\Phi=4+2+1=7,$$
 i.e., for $\hd(E)>7/4$, one has 
$\inter(\Delta_\Phi(E,E,E,E))\ne\emptyset$. 
 \medskip
 
 To verify $(DF)_\sigma$ for $\sigma=(13|24)$, we  calculate
 $$D\Phi=\left[
 \begin{matrix}
 (x^2-x^4)^\perp & (x^4-x^1)^\perp & 0 & (x^1-x^2)^\perp \\
 0 & (x^3-x^4)^\perp & (x^4-x^2)^\perp & (x^2-x^3)^\perp
 \end{matrix}
 \right]$$
 and note that the first and third columns form a matrix of rank two if $x^2\ne x^4$, 
 as do the second and fourth columns under the same condition.
 \medskip

 Pick any $s$ with $7/4<s<\hd(E)$ and let $\mu$ be a Frostman measure on $E$ of finite $s$-energy.
 Arguing as in the earlier proofs, we can
pick a four-tuple $x_0=(x_0^1,x_0^2,x_0^3,x_0^4)$ with each $x_0^j\in supp(\mu)$ such that 

(i) $x_0^2-x_0^4\ne 0$; and

(ii) $x_0^1-x_0^4$ and $x_0^3-x_0^4$ are linearly independent; 

Let $X^j=B\left(x_0^j,\delta\right)$, with $\delta$ chosen small enough so that
(i) and (ii) hold with $x_0$ replaced by any $x\in X^1\times X^2\times X^3\times X^4$. 

Let $\vt_0=(t_0^1,t_0^2)=\Phi(x_0)$. 
Then, we claim that the projections $\pi_L: C^\sigma_{\vt_0}\to T^*X_L$ and $\pi_R: C^\sigma_{\vt_0}\to T^*X_R$ 
drop rank by 1 everywhere; as described in Thm. \ref{thm nondeg fio} (ii),
 it suffices to show this for one of projections, say $\pi_L$. 
 By (ii) above, we can parametrize $Z_{\vt_0}^\sigma$ by $(x^1,x^3,x^4)$, with $x^2$ determined by the 
 nonsingular linear system 
 $$(x^1-x^4)\cdot(x^2-x^4)^{\perp}=t_0^1,\quad (x^2-x^4)\cdot(x^3-x^4)^{\perp}=t_0^2,$$
 whose unique solution we can describe by $x^2=X^2(x^1,x^3,x^4)$. Then
 \bast
 C^\sigma_{\vt_0}&=&\Big\{ \left(x^1,x^3,\theta_1\left(X^2-x^4\right)^\perp, -\theta_2\left(X^2-x^4\right)^{\perp};\,\dots,\dots\right): \\
 & & \qquad\qquad\qquad (x^1,x^3,x^4)\in X^1\times X^3\times X^4, (\theta_1,\theta_2)\in\R^2\setminus 0  \Big\},
 \east
 where the $T^*X^{\hs}$ components on the right are suppressed because they are not needed for the analysis.
 One easily sees that $D\pi_L$ drops rank by 1 everywhere, i.e.,  has constant rank equal to 7, 
 with the image of $\pi_L$ being contained in the hypersurface $\{(x^1,x^3,\xi^1,\xi^3): \xi^1\wedge \xi^3 =0\}$.
 By a fact valid for general canonical relations, $D\pi_R$ also drops rank by $1$ everywhere, as well, and by semicontinuity of the rank,
  $C^\sigma_{\vt}$ drops rank by $k\le 1$ for all $\vt$ close to $\vt_0$. (Thus, $\delta$ above is chosen small enough that all of the values in $\Phi(X)$
  are sufficiently close to $\vt_0$.) 
  By Thm.  \ref{thm nondeg fio}(ii), the $\rvts$ lose at most $\beta_\Phi\le1/2$ derivatives, and we are done.

    
    \subsection{Dot products of differences}\label{subsec prod of diff}
    
To prove Thm. \ref{thm prod of diff},  define 
$$\Phi:\left(\R^d\right)^4\to\R,\, \Phi(x,y,z,w)=(x-y)\cdot(z-w).$$
We will show that using $\sigma=(13|24)$ results in $C_t^\sigma$ which are local canonical graphs, so that
Thm.  \ref{thm kpoint}(ii) applies (with $\beta^\sigma=0$) to yield nonempty interior of the set of dot products of differences for $\hd(E)>(d/2)+(1/4)$.

One computes
$$d\Phi(x,y,z,w)=\left(z-w,\, -\left(z-w\right),\,x-y,\, -\left(x-y\right)\right),$$
so that $(DF)_\sigma$ is satisfied away from $\mathcal W:=\{x-y=z-w=0\}$, which is a codimension $2d$ plane in $\R^{4d}$.
 If $(d/2)+(1/4) <s<\hd(E)$ and $\mu$ is a Frostman measure on $E$ of finite $s$-energy
then, arguing as we have  above,
 $\otimes^4 \mu$ is dominated by $4s$-dimensional Hausdorff measure, and $4s>2d+1$.
Since $\mathcal W$ is a subspace of dimension $2d$,
$\left(\otimes^4\mu \right)(\mathcal W)=0$; repeating previous arguments, we can find base points $x^0,y^0,z^0,w^0\in supp(\mu)$ 
and $\epsilon, \delta>0$ such that $|x-y|+|z-w|>\epsilon$ for $x\in X^1:=B(x^0,\delta),\, y\in X^2:=B(y^0,\delta),\, z\in X^3:=B(z^0,\delta)$ 
and $w\in  X^4:=B(w^0,\delta)$, resp. Thus, $(DF)_\sigma$ is satisfied on $X_L\times X_R$. 
Furthermore, by relabelling and rotating if necessary, we can assume that 
$|z_1-w_1|\ne 0$ on $X_L\times X_R$, so that $d_{x_1}\Phi\ne 0$. 

Thus, letting $t^0=\Phi(x^0,y^0,x^0,w^0)$, for $t$ close to $t^0$, on the hypersurface 
$Z_t^\sigma$  we can solve for $x_1$ as a smooth function of the other variables: $x_1=\bx_1(x',y,z,w)$,
defined for $x'$ in  a small ball $B\subset\R^{d-1}$, and  then parametrize 
\bast
C_t^\sigma&=&\big\{\left(\cdot,\, \cdot,\,\cdot,\,\cdot\,  ; \, y,\, w,\, \theta(z-w),\, \theta( (\bx_1,x')-y)\right)\, \\
& &\qquad\qquad :\, y,z,w\in X^2\times X^3\times X^4,\, x'\in B, \theta\ne 0\big\},
\east
where we have suppressed the $T^*X_L$ entries as irrelevant for the analysis of \linebreak$\pi_R:C_t^\sigma\to T^*X_R=T^*\left(X^2\times X^4\right)$. 
Due to the simple dependence of the  $T^*X^2$ and $X^4$ entries on the coordinates $y,z$ and $w$ on $C^\sigma_t$, 
and denoting elements of  $T^*X^4$ by $(w,\omega)$, we see that
$$rank(D\pi_R)=3d+rank\left(\frac{D\omega}{D(\theta,x')}\right)=4d.$$
Thus, $C^{(13|24)}_t$ is a local canonical graph for $t$ close to $t_0$, and Thm. \ref{thm kpoint} (ii) applies with $k=4,\, p=1,\, d_L=d_R=2d$,
so that for $\hd(E)>(1/4)(2d+1)=(d/2)+(1/4)$, 
$$\inter\left(\left\{(x-y)\cdot(z-w)\, :\, x,y,z,w\in E\, \right\}\right)\ne\emptyset.$$

\subsection{Sum-product sets for bilinear forms}\label{subsec sumprods}

We now state and prove a  more general version of Thm. \ref{thm bilinear forms} on sum-product sets associated to families of bilinear forms.
\ms

\begin{theorem}\label{thm general bilinear forms} Let $\vec{Q}=\left(Q_1,\dots,Q_l\right)$, with
the $Q_j$  nondegenerate, symmetric bilinear forms on $\R^{n_j},\, 1\le j\le l$.
Define $d_1,\dots,d_{2l}$ by $d_{2j-1}=d_{2j}=n_j$,  $1\le j\le l$.
Suppose that   $E_i\subset \R^{d_i}$ are  compact, $1\le i\le 2l$, with 
$$\sum_{i=1}^{2l} \hd(E_i)>1+\frac12\sum_{i=1}^{2l} d_i= 1+\sum_{j=1}^{l} n_j.$$ 
Then the  generalized sum-product set,
\be\label{def form sum set}
\Sigma_{\vec{Q}}\left(E_1,\dots,E_{2l}\right):= 
\left\{ \sum_{j=1}^l Q_j\left(x^{2j-1},x^{2j}\right): x^i\in E_i,\, 1\le i\le 2l\right\}\subset\R,
\ee
has nonempty interior. 
\end{theorem}

Define
$$\Phi\left(x^1,\dots,x^{2l}\right)=\sum_{j=1}^l Q_j\left(x^{2j-1},x^{2j}\right)\hbox{ on }\R^{d_1}\times\cdots\times \R^{d_{2l}}.$$
We show that Thm. \ref{thm general bilinear forms} follows from Thm. \ref{thm kpoint} (ii), 
using $\sigma=\left(13\dots\left(2l-1\right)|24\dots\left(2l\right)\right)$, so that  $d_L=d_R=n:=\sum_{j=1}^ln_j$, and with $p=1$.
Since we may write $Q_j\left(x^{2j-1},x^{2j}\right)=A^jx^{2j-1}\cdot x^{2j}$ for  nonsingular, symmetric $A^j\in R^{n_j\times n_j}$, 
$$d_{x^{2j-1}}\Phi= A^jx^{2j}\hbox{ and } d_{x^{2j}}\Phi= A^jx^{2j-1}.$$
Since the $A^j$ are nonsingular,  all of these are nonzero, 
and thus the double fibration condition \eqref{def DF general} is satisfied if all $x^{2j-1},\, x^{2j}\ne0$. 
Letting $X^i=\R^{d_i}\setminus 0,\, 1\le i\le 2l$, it follows that  $\zts\subset X:=\prod_i X^i$ is a smooth hypersurface,
and we need to analyze the canonical relation in $\left(T^*X_L\setminus 0\right) \times \left( T^*X_R\setminus0\right)$,  
$$C_t^\sigma=\left\{\left(x^1,\,x^3,\dots,x^{2l-3},\,x^{2l-1},\theta A^1x^2,\, \theta A^2x^4,\dots,\theta A^lx^{2l};\dots,\, \dots\right):\, x\in \zts,\, \theta\ne 0\right\},$$
where the entries on the right, in $T^*X_R$, are the even variants of the entries on the left and have been suppressed. 

For each of the $2l$ sets $E_i$, let $s_i<\hd(E_i)$ and $\mu_i$ be a Frostman measure on $E_i$ with finite $s_i$-energy.
Let $E:= E_1\times E_1\times\dots\times E_l\times E_l$ and
pick a base point $x_0:=\left(x^1_0,\dots,\, x^{2l}_0\right)\in E$, 
which we can assume has all of its components nonzero and thus belongs to $X$, 
and a $0<\delta_i< |x^i_0|$ such that  $\mu_i\left(B\left(x_0^i,\delta_i\right)\right)>0$.

Set $t_0=\Phi(x_0)$.
By rotations in $x^1=\left(x^1_1,\dots,x^1_{d_1}\right)=:\left(x^1_1, \left(x^1\right)'\right)$ and $x^2$ if \linebreak necessary, 
we can assume that $d_{x^1_1}\Phi(x_0)\ne 0$,
so that near $x^0$, $Z^\sigma_{t_0}$ is the graph of a function, $x^1_1=f\left(\left(x^1\right)',x^2,\dots,x^{2l}\right)$, with $d_{x^2_1}f\ne 0$.
Hence, we can compute the projection $\pi_L:C^\sigma_{t_0}\to T^*X_L$ with respect to coordinates $\left(x^1\right)',x^2,\dots,x^{2l},\, \theta$.
Since $A^1$ is nonsingular and $\theta\ne 0$, one sees that the map $(x^1)',x^2,\theta$ into the $T^*X^1$ entries has full rank, as do
all of the maps $x^{2j-1},\, x^{2j}$ (with $\theta$ fixed) to $T^*X^{2j-1}$, so that $D\pi_L$ has full rank, and $C^\sigma_{t_0}$ is a local canonical graph.
Hence, $\beta^\sigma=0$  and Thm.  \ref{thm kpoint}(ii) applies, 
yielding $\inter \left(\Sigma_{\vec{Q}}\left(E_1,\dots,E_{2l}\right)\right)\ne\emptyset$ if
$\sum_{i} s_i> n +1$. i.e., if $\sum_i \hd(E_i)>1+(1/2)\sum d_i$.


\section{Final comments}\label{sec comments}

It would be interesting to know whether the Hausdorff dimension thresholds in any of these theorems are sharp.
However, it is worth remarking that the results on pinned volumes and sum-products at least 
have the correct asymptotic behavior as the dimension or the number of quadratic forms tend to infinity,
even for the weaker Falconer problem of positive Lebesgue measure:
\ms

In Thm. \ref{thm volumes}, since all of the volumes are zero if $x^0$ and $E$ both lie in a hyperplane,
one can not take $\hd(E)\le d-1$, and so the restriction $\hd(E)>d-1+(1/d)$ cannot be improved by more than $1/d$.
\ms

Similarly, in Thms. \ref{thm bilinear forms} and \ref{thm general bilinear forms}, 
if we take $E_{2j-1}$ and $E_{2j}$ to be  in $Q_j$-orthogonal subspaces of $\R^d$ 
(in the notation of  Thm. \ref{thm bilinear forms}), 
then $\Sigma_{\vec{Q}}\left(E_1,\dots,E_{2l}\right)=\{0\}$.
Thus,  it is necessary that $\hd(E_{2j-1})+\hd(E_{2j})>d,\, 1\le j\le l$,
so that the $1/l$ in $\hd(E_{2j-1})+\hd(E_{2j})>d+(1/l)$ cannot be reduced  by more than $1/l$.
\ms

Finally, we observe that the results here are obtained by extracting as much as possible from standard estimates for linear Fourier integral operators.
A number of previous results on translation-invariant  Falconer-type configuration problems, such as \cite{GI12,GGIP12,GIM15},
are based on   genuinely  bilinear or multi-linear  estimates for  generalized Radon transforms  and  FIOs, in  settings  where the Fourier transform is an effective tool.
One can ask whether the  thresholds in  this paper (and in \cite{GIT19} for 2-point configurations), 
where the families  $\rvt^{\sigma}$ are  typically nontranslation-invariant,
can be lowered by obtaining truly multi-linear estimates for FIOs.


\medskip

\end{document}